\documentclass[onefignum,onetabnum]{siamart190516}
\usepackage{amsmath,amssymb,amsxtra}
\usepackage{mathrsfs}
\usepackage{float,verbatim}
\usepackage{threeparttable,booktabs}
\usepackage{graphicx}
\usepackage{epstopdf}
\usepackage{mathtools}
\usepackage{subfig}
\usepackage{algorithm, algpseudocode}
\graphicspath{{./pics/}}

\usepackage{array}

\newtheorem{example}{Example}[section]

\usepackage{graphicx}
\usepackage{dsfont}
\usepackage{epstopdf}
\graphicspath{{./pics/}}
\usepackage{booktabs}
\usepackage{threeparttable}
\usepackage{verbatim}
\newcommand{\bel}{\begin{equation} \label}
\newcommand{\ee}{\end{equation}}
\def\beq{\begin{equation}}
\def\eeq{\end{equation}}
\newcommand{\bea}{\begin{eqnarray}}
\newcommand{\eea}{\end{eqnarray}}
\newcommand{\beas}{\begin{eqnarray*}}
\newcommand{\eeas}{\end{eqnarray*}}

\allowdisplaybreaks

\numberwithin{equation}{section}

\allowdisplaybreaks

\def\phi {\varphi}

\allowdisplaybreaks

\allowdisplaybreaks

\title{Iterative Direct Sampling Method for Elliptic Inverse Problems with Limited Cauchy Data\thanks{The work of B. Jin is supported by Hong Kong RGC General Research Fund (14306423) and ANR / RGC Joint Research Scheme (A-CUHK402/24) and a start-up fund from The Chinese University of Hong Kong.
The work of J. Zou was substantially supported by the Hong Kong RGC General Research Fund (projects 14306623 and 14306921) and NSFC/Hong Kong RGC Joint Research Scheme 2022/23 (project N$\_$CUHK465/22).}}

\author{Kazufumi Ito\thanks{Department of Mathematics and Center for Research in Scientific Computation, North Carolina State University, Raleigh, NC 27695-8205 USA. (\texttt{kito@ncsu.edu})} \and
Bangti Jin\thanks{Department of Mathematics, The Chinese University of Hong Kong, Shatin, N.T., Hong Kong (email: \texttt{b.jin@cuhk.edu.hk, fengruwang@cuhk.edu.hk, zou@math.cuhk.edu.hk})}\and
Fengru Wang\footnotemark[3] \and
Jun Zou\footnotemark[3]}

\date{\today}

\begin{document}

\maketitle

\begin{abstract}
In this work, we propose an innovative iterative direct sampling method to solve nonlinear elliptic inverse problems from a limited number of pairs of Cauchy data.
It extends the original direct sampling method (DSM) by incorporating an iterative mechanism, enhancing its performance with a modest increase in computational effort but a clear improvement in its stability against data noise.
The method is formulated in an abstract framework of operator equations and is applicable to a broad range of elliptic inverse problems. Numerical results on electrical impedance tomography, optical tomography and cardiac electrophysiology etc. demonstrate its effectiveness and robustness, especially with an improved accuracy for identifying the locations and geometric shapes of inhomogeneities in the presence of large noise, when compared with the standard DSM.
\end{abstract}



\section{Introduction}
In this work, we are concerned with solving nonlinear inverse problems associated with second-order elliptic partial differential equations (PDEs), and shall develop an innovative iterative direct sampling method (IDSM) that is capable of qualitatively identifying inhomogeneities / inclusions of distributed coefficients in either linear or nonlinear elliptic second-order PDEs using only a limited number of pairs of Cauchy data and enjoys remarkable robustness with respect to data noise.
Elliptic inverse problems of estimating physical parameters are of great importance in various applied fields, e.g., medical imaging (including electrical impedance tomography and diffuse optical tomography etc.), geophysical prospecting, nano-optics, and non-destructive testing / evaluation \cite{Ammari:2008,CheneyNewell:1999,Hecht2006,Schmerr2016,Zhdanov2015}. Since only boundary measurements are available, these inverse problems are frequently severely ill-posed, namely their solutions can be highly susceptible to the presence of data noise, and numerically very challenging to accurately resolve.

Many numerical methods have been proposed to solve these challenging inverse problems. This is often achieved using variational regularization (see, e.g., \cite{Schuster:2012,Jin2014}), i.e., formulating reconstruction tasks as solving suitable optimization problems that incorporate a data-fitting term and a regularization term. The former measures the difference between the model output with the measurement data, whereas the latter combats the inherent ill-posedness. Common regularization terms include Sobolev smoothness \cite{TikhonovArsenin:1977,EnglHankeNeubauer:1996}, sparsity promoting penalty (e.g., $\ell^1$ and elastic net) \cite{Daubechies:2004} and total variation \cite{RudinOsherFatemi:1992}. Numerically, variational regularization leads to a PDE constrained optimization problem. The solutions to the regularized problem are governed by a coupled optimality system, involving the forward PDE, the adjoint PDE, and a variational inequality, and many standard optimization methods then are employed to find the regularized solutions. Alternatively, iterative regularization techniques, e.g., Landweber, Newton, Gauss-Newton or Levernberg-Marquardt methods, can also be employed (together with suitable early stopping rules) \cite{EnglHankeNeubauer:1996}. These methods have been routinely employed in practice with great success, and have been established as standard inversion techniques.

Despite the enormous success in practice, these approaches still suffer from several drawbacks. First, the determination of suitable regularization parameters (and also other algorithmic parameters in the optimizers) is essential for obtaining accurate reconstructions. Yet this remains a very challenging task in general and often requires tedious tuning and searching of these hyper-parameters.
Second, the high nonconvexity of the objective functional may cause the minimization algorithms to get trapped in bad local optima, which necessitates the use of a good initial guess in order to ensure the convergence of the algorithms to the desired global minimizer. However, obtaining a good initial guess can be highly nontrivial.
Note, however, that there has been an  important development on the so-called the convexification method, first proposed by Klibanov in 1997 \cite{Klibanov:1997}. This method can provide globally convergent approximations with explicit convergence estimates, thereby addressing the challenge of local minima and overcoming the drawback of more conventional regularized formulations, and it has been successfully applied to various challenging PDE inverse problems, including inverse acoustic problems, radiative transport equation, and travel time tomography (see, e.g., \cite{Klibanov:1997,Zhipengyang2023,KlibanovLiZhang:2023jcp}); See the recent monograph \cite{KlibanovLi:2021book} for a comprehensive overview.
Third, it can take thousands of iterations for the optimizer to reach convergence.
Additionally, at each step, the gradient of the objective functional with respect to the concerned coefficient(s) is frequently needed to perform the iteration, and the computation of the gradient can be costly for many models involving complex PDEs.

Therefore there is much interest in developing efficient, and robust numerical methods that can provide a rough estimate of inclusions, which may already provide not only valuable information for many practical applications, but can also serve as a reasonable initial guess for iterative methods for further refinements \cite{ItoJinZou:2012jcp}.
Motivated by these observations, direct methods have been intensively investigated during past three decades, which utilize specially designed functionals of the measurement data that can attain high values within inhomogeneities.
These include linear sampling method \cite{Kirsch1996}, multiple signal classification (MUSIC) \cite{Devaney2004,AmmariLesselier:2005,AmmariCalmon:2008}, point source method \cite{Potthast1998}, factorization method \cite{Kirsch1998} and direct sampling method \cite{ItoJinZou:2012} etc.; see the monographs \cite{Monk2011,Potthast2001,Grinberg2008} for recent developments.
Among these methods, the direct sampling method (DSM) has gained much popularity due to its ability to work with a small set of Cauchy data. Initially designed for inverse acoustic and electromagnetic medium scattering \cite{ItoJinZou:2012,ItoJinZou:2013,LiZou:2013}, the DSM has been further developed for scattering problems \cite{HarrisRezac:2020,HarrisNguyen:2020} and extended to non-wave equations, including electrical impedance tomography \cite{Zou2014}, diffuse optical tomography \cite{Zou2015,Zou2021}, Radon transform \cite{ChowHanZou:2021radon}, parabolic type problems \cite{ChowItoZou:2018sisc}, cracks \cite{Park:2018} and Eikonal equation \cite{ItoLiang:2024}, and most recently has been integrated with deep learning techniques for resolution enhancement \cite{GuoJiang:2021,NingHanZou:2024,NingHan:2024b}.
This method employs a carefully chosen Sobolev dual product in the measurement space and a set of probing functions that are nearly orthogonal with respect to this dual product.
By computing the inner product of the difference between the measurement and the solution corresponding to the background medium with the probing functions (often  Green's functions or fundamental solutions), we obtain an index function, which provides qualitative information about the locations and shapes of the inclusions.

Despite its distinct features, the DSM still faces several challenges when applied to elliptic inverse problems.
First, when Green's function is not explicitly available, the DSM necessitates computing the fractional Laplacian of the measurement, which is highly susceptible to data noise.
Second, the method lacks a mechanism to balance accuracy and efficiency, impeding its ability to achieve enhanced results when a higher resolution is needed or when the index function is not sufficiently informative. Third, the DSM indicates the source function induced by the inclusions rather than the inclusions themselves, limiting its accuracy and ability to distinguish inhomogeneities from different physical mechanisms.
To address these technical limitations, we shall develop an iterative direct sampling method (IDSM).
It replaces the Sobolev dual product by a regularized Dirichlet-to-Neumann (DtN) map, thereby markedly enhancing its robustness with respect to noise. It refines the DSM through an iterative approach, enhancing its performance progressively, in order to continuously improve the estimation accuracy.
To maintain the computational efficiency, the IDSM employs a low-rank correction technique, and is scalable to large-scale problems.
Instead of focusing on the induced source, the IDSM directly images the inclusion.
This leads to a significant improvement in the accuracy, and can differentiate inclusions of distinct types.
Moreover, it greatly broadens the scope of its potential application, including complex scenarios that involve non-linearities, e.g., a semilinear model arising in cardiac electrophysiology.
To the best of our knowledge, this is the first direct sampling method for a coefficient inverse problem with a semilinear elliptic model from the Cauchy data. For multiple pairs of Cauchy data, instead of forming an index function for each individual pair like the DSM, the IDSM aggregates the data into one function, which can provide a more informative result that combines the information from all data pairs. In addition, the IDSM is formulated in an abstract framework and is applicable to general elliptic inverse problems. The development of the IDSM in an abstract framework and the extensive numerical illustrations represent the main contributions of the work.

The rest of the paper is organized as follows.
In Section \ref{sec:DSM}, we give an abstract framework for elliptic inverse problems and describe the standard DSM for imaging the inclusions. In Section \ref{sec:IDSM}, we propose the IDSM for elliptic inverse problems in an abstract framework. To illustrate the efficacy and stability of the IDSM, we conduct multiple numerical experiments, including electrical impedance tomography, diffuse optical tomography and electrocardiography, in Section \ref{sec_num}. Throughout, for any Hilbert space $H$ and Banach space $X$, the notation $(\cdot,\cdot)_{H}$ and $\langle \cdot, \cdot\rangle_{X', X}$ denote the inner product of $H$ and the bilinear form on the dual space $X'$ of $X$ and $X$ (i.e., duality product), respectively. For two Sobolev spaces $X$ and $Y$, the notation $B(X,Y)$ denotes the space of all bounded operators that map from $X$ to $Y$, and the set $L(X,Y)\subset B(X,Y)$ consists of all bounded linear operators.

\section{The direct sampling method for elliptic inverse problems}\label{sec:DSM}

In this section, we develop an abstract framework for the concerned class of elliptic inverse problems. The main objective is to identify the locations and shapes of inhomogeneous inclusions within a region using a limited number of pairs of Cauchy data  associated with the elliptic PDE model. To illustrate the idea, we revisit the DSM for electrical impedance tomography (EIT) in Section \ref{subsec_dsm}.
\subsection{Elliptic inverse problems}\label{subsec_eip}

In this section, we formulate an abstract framework for the concerned class of elliptic inverse problems. Throughout, let $\Omega\subset \mathbb{R}^{d}$ ($d = 2, 3$) be an open, bounded, and simply connected domain with a smooth boundary $\Gamma=\partial\Omega$.

First we describe two concrete examples, i.e., electrical impedance tomography (EIT) and the cardiac electrophysiology (CE) problem, which motivate the abstract formulation \eqref{eqn1} below.
\begin{example}\label{exam:EIT}
In EIT, the objective is to estimate the unknown conductivity distribution from voltage measurements taken on the boundary $\Gamma$ {\rm(}see the review \cite{CheneyNewell:1999} for an overview{\rm)}.
Let $\sigma\in L^\infty(\Omega)$ be a strictly positive function representing the conductivity distribution.
The potential $y\in H^{1}(\Omega)$ satisfies the following system:
\begin{equation}\label{eqn26}
\left\{\begin{aligned}
\nabla \cdot(\sigma \nabla y)&=0,\quad \text { in } \Omega, \\
\sigma \frac{\partial y}{\partial n}&=f ,\quad \text { on } \Gamma,
\end{aligned}
\right.
\end{equation}
where $n$ denotes the unit outward normal vector to the boundary $\Gamma$, and $f\in L^{2}(\Gamma)$ represents the current imposed on $\Gamma$, satisfying the compatibility condition $\int_{\Gamma}f\,\mathrm{d}s=0$.
To ensure uniqueness, the solution $y$ is required to belong to the space $\tilde{H}^{1}(\Omega)=\{v\in H^{1}(\Omega):\,\int_{\Gamma}v\,\rm{d}s=0\}$.
The EIT problem aims to determine the regions in which $\sigma$ deviates from a known background conductivity $\sigma_{0}$, using the boundary voltage $y|_{\Gamma}$ {\rm(}often for multiple boundary excitations{\rm)}.
By Green's identity, we can reformulate problem \eqref{eqn26} into a weak form:
\begin{equation}\label{eqn27}
\int_{\Omega}\sigma\nabla y\cdot\nabla v\,{\rm d}x=\int_{\Gamma}fv\,{\rm d}s,\quad \forall v\in\tilde{H}^{1}(\Omega).
\end{equation}
The source $f\in L^2(\Gamma)$ is  the external excitation applied to the boundary $\Gamma$, and can be interpreted as a linear functional on the space $\tilde{H}^{1}(\Omega)$ {\rm(}in view of the trace theorem{\rm)}. The bilinear form $\int_{\Omega}\sigma\nabla y\cdot\nabla v\,{\rm d}x$ on $\tilde{H}^{1}(\Omega)\times \tilde{H}^{1}(\Omega)$ defines a bounded linear operator from $\tilde{H}^{1}(\Omega)$ to $(\tilde{H}^{1}(\Omega))'$.
Let $u=\sigma-\sigma_{0}$ be the inhomogeneity in the conductivity. Then the elliptic operator in \eqref{eqn27} can be decomposed into a background operator that is associated with the background medium and an operator linearly associated with $u$:
\begin{equation*}
\int_{\Omega}\sigma\nabla y\cdot\nabla v\,{\rm d}x=\int_{\Omega}\sigma_{0}\nabla y\cdot\nabla v\,{\rm d}x+\int_{\Omega}u\nabla y\cdot\nabla v\,{\rm d}x,\quad \forall v\in \tilde H^1(\Omega).
\end{equation*}
This decomposition plays a crucial role in developing the DSM and IDSM for EIT below.
\end{example}

\begin{example}\label{exam:EC}
The second example arises in cardiac electrophysiology {\rm(}CE{\rm)} in the context of ischemic heart disease detection \cite{Scacchi2014}.
In this context, the state variable $y$ denotes the electric transmembrane potential, $\omega$ refers to the region occupied by an ischemia, $f$ denotes the current stimulus applied to the tissue and $\widetilde{\sigma}$ is the tissue conductivity, which varies with the presence of ischemia, with $0< \sigma\ll 1$:
\begin{equation*}
\widetilde{\sigma}(x)= \begin{cases}\sigma, & x \in \overline{\omega}, \\ 1, & x \in \Omega \backslash \overline{\omega}.\end{cases}
\end{equation*}
The governing equation of the model reads
\begin{equation}\label{eqn:EC}
\left\{\begin{aligned}
-\nabla \cdot(\widetilde{\sigma} \nabla y)+\chi_{\Omega \backslash \overline{\omega}} y^3 & =f, & & \text { in } \Omega, \\
\widetilde{\sigma}\frac{\partial y}{\partial n} & =0, & & \text { on } \Gamma,
\end{aligned}\right.
\end{equation}
where $\chi_{\Omega \backslash \overline{\omega}}$ is the characteristic function of the set $\Omega \backslash \overline{\omega}$.
The cubic nonlinearity in the model \eqref{eqn:EC} represents one significant empirical law.
In cardiac electrophysiology, one aims at locating the ischemia region $\omega$ using the electric voltage measurement $y|_{\Gamma}$ on the boundary $\Gamma$.
By Green's identity, we can recast problem \eqref{eqn:EC} into a weak form:
\begin{equation*}
\int_{\Omega}\widetilde{\sigma}\nabla y\cdot\nabla v\,{\rm d}x+\int_{\Omega}\chi_{\Omega\backslash\bar\omega}y^{3}v\,{\rm d}x=\int_{\Omega}fv\,{\rm d}x,\quad \forall v\in H^{1}(\Omega).
\end{equation*}
Like before, $f$ is a linear functional on the space $H^{1}(\Omega)$, but $ \int_{\Omega}\widetilde{\sigma}\nabla y\cdot\nabla v\,{\rm d}x+\int_{\Omega}\chi_{\Omega\backslash\bar\omega}y^{3}v\,{\rm d}x$ defines a bounded semilinear elliptic operator from $H^{1}(\Omega)$ to $\left(H^{1}(\Omega)\right)'$.
Upon letting $u=\chi_{\overline\omega}$, the semilinear elliptic operator can also be decomposed into an operator corresponding to the background and an operator linearly associated with $u$:
\begin{align}
{}&\int_{\Omega}\widetilde{\sigma}\nabla y\cdot\nabla v\,{\rm d}x+\int_{\Omega}\chi_{\Omega\backslash\overline\omega}y^{3}v\,{\rm d}x\nonumber\\
=&\int_{\Omega}\left[1+u(\sigma-1)\right]\nabla y\cdot\nabla v\,{\rm d}x+\int_{\Omega}(1-u)y^{3}v\,{\rm d}x\nonumber\\
=&\int_{\Omega}\nabla y\cdot\nabla v + y^{3}v\,{\rm d}x+\int_{\Omega}u\left[(\sigma - 1)\nabla y\cdot \nabla v - y^{3}v\right]\,{\rm d}x.\label{eqn16}
\end{align}
This relation forms the basis for developing the IDSM for the model.
\end{example}

These two examples motivate formulating an abstract framework for the concerned class of elliptic inverse problems, in which the objective is to determine the locations and shapes of the inhomogeneities within the background medium. We denote by $u$ one or multiple physical parameters that generate inhomogeneous inclusions, characterized by their distinct material properties compared to the background. Typically, the function $u\in S\subset L^{\infty}(\Omega)$ for one single type of inhomogeneity may take the following form
\begin{equation*}
u(x) = \sum_{i=1}^nc_i \chi_{\omega_i}(x),
\end{equation*}
with $n$ being the number of inclusions, the subdomains $\omega_i\subset \Omega$ mutually disjoint from each other (and also away from the boundary $\Gamma$), and $c_i\in \mathbb{R}$ a (possibly unknown) constant coefficient(s) within $\omega_i$. The notation $\chi_{\omega_i}$ denotes the characteristic function of the set $\omega_i$. Often the admissible set $S$ encodes various physical constraints. For example, the conductivity $\sigma$ in EIT is typically strictly positive with a positive lower bound.  Then consider the following direct problem for the following (possibly nonlinear) second-order elliptic equation in an operator form:
\begin{equation}\label{eqn1}
\mathcal{A}[y]y + \mathcal{B}[u](y) = f,
\end{equation}
where $f\in \left(H^{1}(\Omega)\right)'$ represents either a boundary or domain source.
The elliptic operator $\mathcal{A}:H^{1}(\Omega)\to B(H^{1}(\Omega),\left(H^{1}(\Omega)\right)')$ models the physical process in the background (i.e., in the absence of the inclusions) and can be nonlinear in the state $y$, whereas the operator $\mathcal{B}$ describes how the inhomogeneity $u$ enters into the problem formulation and influences the state variable $y$, and linearly maps from $S\subset L^{\infty}(\Omega)$ to $B(H^{1}(\Omega), \left(H^{1}(\Omega)\right)')$. This decomposition isolates the inhomogeneity $u$ on the state $y$, which facilitates developing direct methods for imaging $u$ from the boundary measurements.

Note that the operator equation \eqref{eqn1} should be understood in a distributional sense. For Examples \ref{exam:EIT} (EIT) and \ref{exam:EC} (CE), the operators $\mathcal{A}$ and $\mathcal{B}$ are given respectively by
\begin{align*}
\langle \mathcal{A}[y](y),v\rangle_{(H^1(\Omega))',H^1(\Omega)} &  = \left\{\begin{aligned}
\int_{\Omega}\sigma_{0}\nabla y\cdot\nabla v\,{\rm d}x, &\quad \mbox{EIT},\\
\int_{\Omega}\nabla y\cdot\nabla v + y^{3}v\,{\rm d}x, & \quad \mbox{CE},
\end{aligned}\right.\\
\langle \mathcal{B}[u](y),v\rangle_{(H^1(\Omega))',H^1(\Omega)} &  = \left\{\begin{aligned}
\int_{\Omega}u\nabla y\cdot\nabla v\,{\rm d}x, &\quad \text{EIT},\\
\int_{\Omega}u\left[(\sigma - 1)\nabla y\cdot \nabla v - y^{3}v\right]\,{\rm d}x,&\quad\mbox{CE},
\end{aligned}\right.
\end{align*}
where $v\in \tilde H^1(\Omega) / H^1(\Omega) $ is a test function.
Thus, the abstract formulation \eqref{eqn1} covers these two examples. More examples will be given in Section \ref{sec_num}.

The concerned inverse problem reads: to recover the shapes and locations of inhomogeneous inclusions $\omega_i$ given one or a small number of pairs of Cauchy data $(f, y_{d})$, where $y_{d}$ is the Dirichlet data of the solution $y$ to problem \eqref{eqn1} on the boundary $\Gamma$.
This task represents a fundamental challenge in practice due to the severe ill-posed nature. We shall develop a novel iterative DSM for efficiently imaging the inclusions $\omega_i$.

\subsection{The direct sampling method}\label{subsec_dsm}
Now we describe the general principle of the DSM in the abstract framework, and illustrate it on EIT in Example \ref{exam:EIT} with one pair of Cauchy data. For multiple pairs of Cauchy data, one may compute the index function individually for each pair, and then form their sum as the overall index.
Throughout, $\mathcal{T}\in L(H^{1}(\Omega), L^{2}(\Gamma))$ denotes the trace operator, where the choice of the space $L^2(\Gamma)$ (instead of the natural choice $H^{1/2}(\Gamma)$) as the range of $\mathcal{T}$ is to match that of the noisy boundary data.
We denote by $u_{*}$ and $y_{*}$ the ground truth for $u$ and $y$, and by $y_{d}$ the noisy measurement, which is a corrupted version of $\mathcal{T}y_{*}$.
Recall that in the context of EIT, the inhomogeneity $u_{*}$ and the ground-truth potential $y_{*}$ satisfy the following system
\begin{equation*}
\left\{\begin{aligned}
-\nabla \cdot\left[\left(\sigma_{0}+u_{*}\right) \nabla y_{*}\right]&=0, \quad \text {in } \Omega, \\
\sigma_{0} \frac{\partial }{\partial n}y_{*}&=f, \quad\text {on } \Gamma,
\end{aligned}
\right.
\end{equation*}
if the inclusion $u_*$ is compactly supported in $\Omega$ (i.e., $\overline{\text{supp}(u)}\cap\Gamma=\emptyset$).
We denote the potential for the homogeneous model (corresponding to the background $\sigma_0$) by $y_{\emptyset}$, which solves
\begin{equation*}
\left\{\begin{aligned}
-\nabla \cdot\left(\sigma_{0} \nabla y_{\emptyset}\right)&=0, \quad \text {in } \Omega, \\
\sigma_{0} \frac{\partial }{\partial n}y_{\emptyset}&=f, \quad\text {on } \Gamma.
\end{aligned}
\right.
\end{equation*}
The difference between the potentials of these two models on the boundary $\Gamma$, i.e.,
\begin{equation*}
y^{s}_* := \mathcal{T}(y_{\emptyset}-y_{*}),
\end{equation*}
is called the scattering field below.
With the induced source $F=\nabla\cdot(u_{*}\nabla y_{*})$, one has the following equation for $y_{\emptyset}-y_{*}$:
\begin{equation}\label{eqn14}
\left\{\begin{aligned}
-\nabla \cdot\left[\sigma_{0} \nabla \left(y_{*}-y_{\emptyset}\right)\right]&=F, \quad \text {in } \Omega, \\
\sigma_{0} \frac{\partial }{\partial n}\left(y_{*}-y_{\emptyset}\right)&=0, \quad\text {on } \Gamma.
\end{aligned}
\right.
\end{equation}
The difference $y_{\emptyset}-y_{*}$ is linearly related to the source $F$ via an integral representation. We denote Green's function (Neumann's function, more precisely) of the homogeneous model by $G(x, x')$:
\begin{equation*}
\left\{\begin{aligned}
{}-\nabla_{x'}\cdot(\sigma_{0}\nabla_{x'} G(x,x'))&=\delta(x'-x),\quad \text{in }\Omega,\\
\sigma_0\frac{\partial{}}{\partial n_{x'}}G(x,{x'}) &= -\frac{1}{|\Gamma|},\quad \text{on }\Gamma,\\
\int_{\Gamma}G(x,x')\,\mathrm{d}s_{x'}&=0,
\end{aligned}\right.
\end{equation*}
where the notation ${\rm d}s_{x'}$ denotes the infinitesimal arc length in the variable $x'$, and $\delta(x)$ denotes the Dirac delta function.
Formally $G(x,x')$ is the response of the system at $x'$ to a point source located at $x$. It can be used to represent the scattering field $y_*^s$ as
\begin{equation}\label{eqn:scattering}
y_*^s(x)=\int_{ \Omega}G(x,x^{\prime})F(x^{\prime})\,\mathrm{d}x^{\prime}\approx \sum_{i}G(x,x_{i})F(x_{i}),\quad \text{for }x\in \Gamma.
\end{equation}
That is, the scattering field $y_*^s$ can be well approximated by a linear combination of $G(x,x')$, with the weight given by the source $F$. Namely, $G(x,x')$ is fundamental in the sense that it can well represent the scattering field $y_*^s$. By incorporating the inclusions $u$ into the source $F$, one can identify the locations $u$ by analyzing the distribution of $F$.
This fact is foundational to the DSM when identifying inclusions.
By pinpointing the locations at which the source $F$ is significant, the locations of the inclusions can be inferred.
To this end, we need a dual product $\langle \cdot, \cdot\rangle_{\Gamma}$ on the data space. This dual product $\langle \cdot, \cdot\rangle_{\Gamma}$ is designed to ensure that Green's functions $G(\cdot,x_i)$ exhibit the property of nearly mutual orthogonality when their arguments $x_i$s are well separated:
\begin{equation*}
\frac{\langle G(\cdot, x_{1}), G(\cdot, x_{2})\rangle_{\Gamma}}{\langle G(\cdot, x_{1}), G(\cdot, x_{1})\rangle_{\Gamma}}\ll 1,\quad\text{ if }|x_{1}-x_{2}|\text{ is large.}
\end{equation*}
This property is crucial for the success of the following DSM index function $\eta$, defined by
\begin{equation}\label{eqn2}
\eta(x)=\frac{\langle G(\cdot, x), y_d^s \rangle_{\Gamma}}{\langle G(\cdot, x), G(\cdot, x)\rangle_{\Gamma}}\approx\frac{\langle G(\cdot, x), y_*^s\rangle_{\Gamma}}{\langle G(\cdot, x), G(\cdot, x)\rangle_{\Gamma}}\approx \sum_{i}\frac{\langle G(\cdot, x), G(\cdot, x_{i})\rangle_{\Gamma}}{\langle G(\cdot, x), G(\cdot, x)\rangle_{\Gamma}}F(x_{i}),
\end{equation}
where we use the noisy scattering field $y_d^s$ defined by
\begin{equation*}
    y_d^s : = \mathcal{T}y_\emptyset - y_d,
\end{equation*}
upon replacing the exact data $\mathcal{T}y_*$ with the noisy measurement $y_{d}$.
The function $\eta(x)$ indicates the presence of $F$: it takes value of order $O(1)$ when the sampling point $x$ is in the vicinity of the set $\{x_{i}\}\subset\mathrm{supp}(F)$ and its value is small otherwise.

Finding a suitable dual product $\langle \cdot, \cdot\rangle_{\Gamma}$ for $G(x,x')$ is a challenging task.
For wave-type problems, the standard $L^{2}(\Gamma)$ product is commonly utilized, and works fairly well \cite{ItoJinZou:2012,ItoJinZou:2013,LiZou:2013}.
However, it does not work so well for elliptic problems, as it does not satisfy the property of nearly mutual orthogonality.
To overcome this limitation, in a series of important works \cite{Zou2014,Zou2015,Zou2021,Zou2022}, Chow et al. propose to use the $H^{\gamma}(\Gamma)$ inner product.
Specifically, let $\{(\lambda_{i},\phi_{i})\}_{i=1}^{\infty}$ be the eigenpairs of Laplace-Beltrami operator $-\Delta_{\Gamma}$ in $L^{2}(\Gamma)$ with $\{\phi_i\}_{i=1}^\infty$ forming an orthonormal basis in $L^2(\Gamma)$.
Then the action of the operator $(-\Delta_\Gamma)^s$ and $H^{\gamma}(\Gamma)$ inner product $(\cdot, \cdot)_{H^{\gamma}(\Gamma)}$ are given by
\begin{equation*}
(-\Delta_\Gamma)^{\gamma}y=\sum_{i=1}^{\infty}\lambda_{i}^{\gamma}(y,\phi_i)_{L^{2}(\Gamma)}\phi_{i}\quad \mbox{and}\quad( y_{1},y_{2})_{H^{\gamma}(\Gamma)}=\int_{\Gamma}\left[(-\Delta_\Gamma)^{\gamma/2}y_{1}\right]\left[(-\Delta_\Gamma)^{\gamma/2}y_{2}\right]\,\mathrm{d}s,
\end{equation*}
respectively. Note that Green's functions $G(x,x')$ are generally smooth on the boundary $\Gamma$ for smooth domains (if $\sigma_{0}$ is also smooth). Therefore even when the noisy scattering field
$y_d^s$ is not smooth, the use of the $H^{\gamma}(\Gamma)$ product is still feasible, since by integration by parts, one can transfer differentiation from the function $y_d^s$ to $G(x,x')$ (the subscript $x'$ indicates that the operation is in the variable $x'$):
\begin{equation*}
\begin{aligned}
(G(x,\cdot), y_d^s)_{H^{\gamma}(\Gamma)}=&\int_{\Gamma}\left[(-\Delta_{\Gamma,x'})^{\gamma/2}G(x,x')\right]\left[(-\Delta_{\Gamma,x'} )^{\gamma/2}y_d^s(x')\right]\,\mathrm{d}s_{x'}\\
=&\int_{\Gamma}\left[(-\Delta_{\Gamma,x'})^{\gamma}G(x,x')\right]y_d^s(x')\,\mathrm{d}s_{x'}\\
=&( (-\Delta_{\Gamma,x'})^{\gamma}G(x,x'),y_d^s)_{L^{2}(\Gamma)}.
\end{aligned}
\end{equation*}
This fact is essential to achieve the robustness of the DSM in the presence of large data noise.

The other critical aspect is the evaluation of the index function $\eta$.
When the domain $\Omega$ is simple, e.g., a disk or a rectangle, $G(x,x')$ can be expressed explicitly using the method of separation of variables.
However, when the explicit form of $G(x,x')$ is unavailable, it is not necessary to solve one direct problem for each sampling point; instead, one single solution suffices.
To see this, let $\zeta\in \tilde{H}^{1}(\Omega)$ be the solution of the following problem:
\begin{equation}\label{eqn3}
\left\{\begin{aligned}
-\nabla\cdot(\sigma_{0}\nabla\zeta)&=0, & \text { in } \Omega, \\
\sigma_0 \frac{\partial \zeta}{\partial n}&=(-{\Delta_{\Gamma}})^{\gamma}y^{s}_{d}, & \text { on } \Gamma.
\end{aligned}
\right.
\end{equation}
It follows from Green's identity that
\begin{align*}
&\int_{\Omega}G(x,x')\nabla\cdot(\sigma_0\nabla\zeta(x)) -\zeta(x)\nabla\cdot(\sigma_0\nabla G(x,x'))\,\mathrm{d}x\\
=&\int_{\Gamma}G(x,x')\sigma_0\frac{\partial}{\partial n_x}\zeta(x)-\zeta(x)\sigma_0\frac{\partial}{\partial n_x}G(x,x') \,\mathrm{d}s_{x}.
\end{align*}
This, along with the definition of $G(x,x')$ and \eqref{eqn3}, implies that $\zeta(x)$ can be expressed as
\begin{align*}
\zeta(x)=&\int_{\Gamma}G(x,x^{\prime}){\left[(-{\Delta_{\Gamma,x'}})^{\gamma}y^{s}_{d}(x^{\prime})\right]}\,\mathrm{d}s_{x'}\\
=&\int_{\Gamma}{\left[(-{\Delta_{\Gamma,x'}})^{\gamma/2}G(x,x^{\prime})\right]\left[(-{\Delta_{\Gamma,x'}})^{\gamma/2}y^{s}_{d}(x^{\prime})\right]}\,\mathrm{d}s_{x^{\prime}}\\
=&( G(\cdot, x), y^{s}_{d})_{H^{\gamma}(\Gamma)}.
\end{align*}
This strategy only requires solving one direct problem to determine the numerator of $\eta$. For the denominator, an approximation based on the distance function $d(x, \Gamma):=\inf_{x'\in\Gamma}|x-x'|$ (with $|\cdot|$ being the Euclidean norm) is often employed:
\begin{equation}\label{eqn22}
\langle G(\cdot, x), G(\cdot, x)\rangle_{\Gamma}\approx Cd(x, \Gamma)^{\gamma},
\end{equation}
where the exponent $\gamma$ is often determined empirically.
Since the inclusion is indicated by the region in which $\eta$ takes relatively large values, only the relative magnitude of $\eta$ is needed (and thus the constant $C$ may be omitted).
Alternatively, one may use the fundamental solution $\Phi_x$ of the Laplace equation $-\Delta\Phi_{x}=\delta_{x}$ in $\mathbb{R}^{d}$, i.e., $\Phi_{x}(x^{\prime})=-\frac{1}{2\pi}\ln|x-x^{\prime}|$ in $\mathbb{R}^{2}$ and $\Phi_{x}(x^{\prime})=\frac{1}{4\pi}|x-x^{\prime}|^{-1}$ in $\mathbb{R}^{3}$.
Then the norm of $\Phi_{x}$ on the boundary $\Gamma$ can also approximate the denominator:
\begin{equation*}
( G(\cdot, x), G(\cdot, x))_{H^\gamma(\Gamma)}\approx \|\Phi_{x}\|_{H^{\gamma}(\Gamma)}^{2},
\end{equation*}
due to its similar asymptotic behavior to that of Green's function:
$$\lim_{x{\to} x'}\frac{G(x,x')}{\Phi_{x}(x')}=\sigma_0(x')^{-1}.$$

\section{Iterative direct sampling method}\label{sec:IDSM}
Despite its salient efficiency and stability, the DSM still has several limitations.
The first issue is its limited accuracy.
The potential causes include lack of near orthogonality, i.e., the inner products between $G(\cdot,x')$s do not decay sufficiently fast, or an inaccurate estimation of $\|G(x,\cdot)\|_{H^\gamma(\Gamma)}$.
Furthermore, the DSM lacks a mechanism to enhance its accuracy.
The second issue is that $\eta$ is designed to image the source $F$ rather than the inclusion $u_*$ directly.
For example, in EIT, the source $F$ is given by $F=\nabla\cdot(u_{*}\nabla y_{*})$, and if $u_*$ is the characteristic function $\chi_\omega$ of a region $\omega$, $F$ contains a singularity supported on the inclusion boundary $\partial\omega$, which can diminish the distinguishability of $\eta$.
Furthermore, the source $F$ may originate from  inclusions of distinct types \cite{Zou2022} (see also Example 2 in Section \ref{sec_num}).
The DSM does not allow differentiating them within the support of $F$.
The third issue is noise robustness.
When an explicit form of $G(x,x')$ is unavailable, problem \eqref{eqn3} has to be solved, which involves computing the fractional Laplacian $(-\Delta_{\Gamma})^\gamma y^{s}_{d}$ of the noisy scattering field $y^{s}_{d}$.
This step essentially amounts to numerical differentiation, and may lead to significant errors, especially in the presence of large random noise.

To overcome these limitations of the DSM, we shall develop a novel iterative direct sampling method (IDSM).
Formally, it aims to
continuously refine $\eta$ with only a modest increase in the computational cost by iterating the DSM process. This ambitious goal is achieved via the following three key algorithmic innovations:
\begin{itemize}
\item[(i)] To replace the dual product $(-\Delta_\Gamma)^\gamma$ by the regularized Dirichlet to Neumann map (DtN map) $\Lambda_\alpha(\mathcal{A})$; see \eqref{eqn7} below for the definition.
It effectively filters out the noise in the data $y^{s}_{d}$, and improves the robustness.
\item[(ii)] To integrate the operator $\mathcal{B}$ into the forward operator.
It redirects $\eta$ to focus on the recovery of the inhomogeneity $u$ directly instead of the induced source $F$.
\item[(iii)] To adjust $\eta$ using a low-rank correction technique.
This not only greatly improves the accuracy of $\eta$ but also maintains the computational efficiency.
\end{itemize}

Note that the first two refinements complicate the forward equation, cf. \eqref{eqn14}, making the derivation of the associated Green's function $G(x,x')$ or its approximations more challenging. The third innovation is to address the challenge.
In addition, the {iterative scheme} also enables the method to adjust the denominator of $\eta$ in \eqref{eqn2} based on a regularized DtN map $\Lambda_\alpha(\mathcal{A})$, allowing for a broad range of hyperparameter values, thereby enhancing its flexibility.

The rest of the section is organized as follows. In Section \ref{ssec:reg-DtN}, we describe the main technical tools, i.e., the regularized DtN map $\Lambda_\alpha(\mathcal{A})$ and the operator $B_\tau$. Then we recast the DSM in an operator form in Section \ref{ssec:DSM}, and finally derive the IDSM in Section \ref{ssec:IDSM} that implements these three innovations.

\subsection{Regularized DtN map}\label{ssec:reg-DtN}
We first introduce some useful tools to construct the IDSM. In the IDSM, the regularized DtN map replaces the fractional Laplacian $(-\Delta_{\Gamma})^\gamma$ on the boundary $\Gamma$ (see \eqref{eqn3}).
For any $v\in L^{2}(\Gamma)$, the DtN map, also known as the Poincar\'e–Steklov operator, is given by $\Lambda_{0}v = \frac{\partial}{\partial n}w$, where $w$ solves
\begin{equation}\label{eqn4}
\left\{
\begin{aligned}
-\Delta w&=0,&\text{ in }\Omega,\\
w&=v, &\text{ on }\Gamma.
\end{aligned}
\right.
\end{equation}
The DtN map $\Lambda_{0}$ resembles the square root {$(-\Delta_{\Gamma})^{1/2}$} of the Laplace-Beltrami operator $-\Delta_{\Gamma}$, exhibiting similar spectral behavior and principal symbols \cite{Polterovich2021}.
It can thus be used as an alternative to $(-\Delta_{\Gamma})^\gamma$  ($\gamma=\frac{1}{2}$).
However, the spectral behavior of the DtN map $\Lambda_0$ implies that it is susceptible to $L^2(\Gamma)$ noise: as an unbounded operator on $L^2(\Gamma)$, its action on the noisy scattering field $y^{s}_{d}$ can amplify the noise greatly.
To remedy this issue, we regularize $\Lambda_0$ so that it is stable in $L^{2}(\Gamma)$, but also can approximate the DtN map $\Lambda_{0}$, hence also $(-\Delta_{\Gamma})^{1/2}$.

Using Green's identity, the Ritz variational formulation of \eqref{eqn4} reads:
\begin{equation}\label{eqn:Ritz}
w\in \mathop{\arg\min}_{w\in H^{1}(\Omega)}\frac12\|\nabla w\|^{2}_{L^{2}(\Omega)}\quad \text{subject to }\mathcal{T}w=v.
\end{equation}
For an arbitrary $v\in L^2(\Gamma)$, problem \eqref{eqn:Ritz} may not have a solution $w\in H^1(\Omega)$.
Formally, the Lagrangian $\mathcal{L}(w,z)$ of problem \eqref{eqn:Ritz} is given by
\begin{equation*}
\mathcal{L}(w,p)=\frac{1}{2}\|\nabla w\|^{2}_{L^{2}(\Omega)}-\left(p,\mathcal{T}w-v\right)_{L^{2}(\Gamma)},
\end{equation*}
with $p\in L^{2}(\Gamma)$ being the Lagrange multiplier. The Karush-Kuhn-Tucker conditions are given by
\begin{equation}\label{eqn18}
\left\{
\begin{aligned}
(\nabla w, \nabla z)_{L^{2}(\Omega)^{d}}-(p, \mathcal{T}z)_{L^{2}(\Gamma)}&=0,& \forall z\in H^{1}(\Omega),\\
(\mathcal{T}w,q)_{L^{2}(\Gamma)}&=(v,q)_{L^{2}(\Gamma)},& \forall q\in L^{2}(\Gamma).
\end{aligned}
\right.
\end{equation}

It follows from the first row of \eqref{eqn18} that $p = \frac{\partial}{\partial n}w$, i.e., the output of the DtN map $\Lambda_{0}$.
To overcome the instability of the system, we employ the standard regularization technique and add a penalty term $\|p\|_{L^2(\Gamma)}^{2}$ (with the penalty parameter $\alpha>0$) to the Lagrangian $L(w, p)$:
\begin{equation*}
\mathcal{L}_{\alpha}(w,p)=\frac{1}{2}\|\nabla w\|^{2}_{L^{2}(\Omega)}-\left(\mathcal{T}w-v,p\right)_{L^{2}(\Gamma)}+\frac{\alpha}{2} \|p\|_{L^{2}(\Gamma)}^{2}.
\end{equation*}
Then the KKT conditions read
\begin{equation}\label{eqn6}
\left\{\begin{aligned}
(  \nabla w, \nabla z)_{L^{2}(\Omega)^{d}}-( p,\mathcal{T}z)_{L^{2}(\Gamma)}&=0,& \forall z\in H^{1}(\Omega),\\
(\mathcal{T}w,q)_{L^{2}(\Gamma)}+(\alpha p,q)_{L^{2}(\Gamma)}&=(v,q)_{L^{2}(\Gamma)},& \forall q\in L^{2}(\Gamma).
\end{aligned}\right.
\end{equation}
We denote by $\Lambda_{\alpha}v=p$ the solution operator of this system.
By the standard regularization theory, the operator $\Lambda_{\alpha}$ approximates the DtN map for small $\alpha>0$, and is stable in the $L^{2}(\Gamma)$ norm.
By Green’s identity, we can see that $w$ solves the following Robin boundary value problem:
\begin{equation*}
\left\{
\begin{aligned}
-\Delta w&=0,&\text{ in }\Omega,\\
\frac{\partial}{\partial n}w-p&=0,&\text{ on }\Gamma,\\
w+\alpha p&=v,&\text{ on }\Gamma.
\end{aligned}
\right.
\end{equation*}
The regularized DtN map $\Lambda_{\alpha}$ is weakly continuous with respect to $\alpha$ \cite{Wang2024}:
\begin{equation*}
\lim_{\alpha \to 0^+}\left( \left(\Lambda_{\alpha}-\Lambda_{0}\right)v,v\right)_{L^{2}(\Gamma)}=0,\quad\forall v\in C^{\infty}(\Gamma).
\end{equation*}
Thus it can approximate the DtN map $\Lambda_{0}$ as well as $(-\Delta_{\Gamma})^{1/2}$.
Furthermore, $\Lambda_{\alpha}$ enjoys excellent stability in $L^2(\Gamma)$ (see Lemma \ref{lemma1} below for the precise statement), and thus is amenable with the presence of $L^2(\Gamma)$ noise, which contrasts sharply with $(-\Delta_{\Gamma})^{1/2}$.
This is the main motivation to adopt $\Lambda_\alpha$ (with a small $\alpha$ value) to approximate $(-\Delta_{\Gamma})^{1/2}$.

In practical implementation, employing \eqref{eqn6} can lead to increased computational complexity.
Indeed, when computing the indicator function $\eta$, \eqref{eqn6} is solved twice for both the forward operator and the dual operator, as illustrated in \eqref{eqn10}.
However, Lemma~\ref{lemma2} below shows that the computational overhead can be mitigated by replacing the kernel operator $\nabla^{*}\nabla$ with the elliptic operator $\mathcal{A}$ of the model \eqref{eqn1}, which requires solving equation \eqref{eqn7} only once.
Thus, in the IDSM, we use the regularized DtN map $\Lambda_\alpha{(\mathcal{A})}$. For any $v\in L^{2}(\Gamma)$, $\alpha> 0$ and the linear elliptic operator $\mathcal{A}\in L(H^{1}(\Omega),\left(H^{1}(\Omega)\right)')$, we define the regularized DtN map $\Lambda_{\alpha}(\mathcal{A})$ to be $\Lambda_{\alpha}(\mathcal{A})v=p$, where the tuple $(w,p)\in H^{1}(\Omega)\times L^{2}(\Gamma)$ solves
\begin{equation}\label{eqn7}
\left\{
\begin{aligned}
\langle \mathcal{A}w, z\rangle_{(H^{1}(\Omega))', H^{1}(\Omega)}-( p,\mathcal{T}z)_{L^{2}(\Gamma)}&=0,& \forall z\in H^{1}(\Omega),\\
(\mathcal{T}w,q)_{L^{2}(\Gamma)}+(\alpha p,q)_{L^{2}(\Gamma)}&=(v,q)_{L^{2}(\Gamma)},& \forall q\in L^{2}(\Gamma).
\end{aligned}
\right.
\end{equation}

The next lemma shows that the regularized DtN map $\Lambda_{\alpha}(\mathcal{A})$ is well defined, with an explicit stability estimate in terms of the Dirichlet boundary data.
\begin{lemma}\label{lemma1}
Fix $\alpha>0$ and let $\mathcal{A}\in L(H^{1}(\Omega),\left(H^{1}(\Omega)\right)')$ be a positive definite operator in the sense that there exist $m, M>0$ such that
\begin{align*}
m\|w\|_{H^{1}(\Omega)}^{2}\leq&\langle \mathcal{A}w, w\rangle_{(H^{1}(\Omega))', H^{1}(\Omega)},\quad \forall w\in H^1(\Omega),\\
\langle \mathcal{A}w,z\rangle_{(H^{1}(\Omega))', H^{1}(\Omega)}\leq& M\|w\|_{H^{1}(\Omega)}\|z\|_{H^{1}(\Omega)},\quad \forall w,z\in H^1(\Omega).
\end{align*}
Then there exists a unique solution $p\in L^2(\Omega)$ to \eqref{eqn7} and moreover, $\|p\|_{L^{2}(\Gamma)}\leq \frac{C}{\min(\alpha,m)} \|v\|_{L^{2}(\Gamma)}$.
\end{lemma}
\begin{proof}
We denote the bilinear form and linear functional associated with \eqref{eqn7} by $a([w,p],[z,q])$ and $f([z,q])$, respectively, given by
\begin{align*}
a([w,p],[z,q])&=\langle \mathcal{A}w, z\rangle_{(H^{1}(\Omega))', H^{1}(\Omega)}-( p,\mathcal{T}z)_{L^{2}(\Gamma)}+( \mathcal{T}w,q)_{L^{2}(\Gamma)}+( \alpha p,q)_{L^{2}(\Gamma)},\\
f([z,q])&=(v,q)_{L^{2}(\Gamma)},
\end{align*}
for all $[w,p],[z,q]\in H^{1}(\Omega)\times L^{2}(\Gamma)$.
The continuity of the elliptic operator $\mathcal{A}$ and the trace operator $\mathcal{T}$ yield the continuity of $a([w,p],[z,q])$ and $f([z,q])$.
The coercivity of $a([w,p],[z,q])$ follows by
\begin{align*}
a([w,p],[w,p])=&\langle \mathcal{A}w, w\rangle_{(H^{1}(\Omega))', H^{1}(\Omega)}-( p,\mathcal{T}w)_{L^{2}(\Gamma)}+( \mathcal{T}w,p)_{L^{2}(\Gamma)}+( \alpha p,p)_{L^{2}(\Gamma)}\\
=&\langle \mathcal{A}w,w\rangle_{(H^{1}(\Omega))', H^{1}(\Omega)}+( \alpha p,p)_{L^{2}(\Gamma)}\\
\geq &m\|w\|_{H^{1}(\Omega)}^{2}+\alpha \|p\|_{L^{2}(\Gamma)}^{2}\geq \min(m,\alpha)(\|w\|_{H^1(\Omega)}^2+\|p\|_{L^2(\Gamma)}^2).
\end{align*}
This and the Lax-Milgram theorem directly give the well-posedness of problem \eqref{eqn7} and the stability estimate $$\|p\|_{L^2(\Gamma)}\leq \frac{C}{\min(\alpha,m)}\|v\|_{L^2(\Gamma)},$$ completing the proof of the lemma.
\end{proof}

Moreover, we define the operator $\mathcal{B}_{\tau}$ to be the adjoint of $\mathcal{B}$ as
\begin{equation}\label{eqn17}
\begin{aligned}
\mathcal{B}_{\tau}:&H^{1}(\Omega){\to} L\left(L^{\infty}(\Omega), \left(H^{1}(\Omega)\right)'\right),\\
{}&\mathcal{B}_{\tau}[y]u=\mathcal{B}[u](y).
\end{aligned}
\end{equation}
The operator $\mathcal{B}_\tau$ allows expressing the source $F$ as $F = \mathcal{B}_{\tau}[y]u$. Note that the operator $\mathcal{B}_\tau$ is the transpose of $\mathcal{B}$ in the arguments $y$ and $u$, and the subscript $\tau$ is to distinguish it from the adjoint operator $\mathcal{B}_{\tau}[y]^{*}$ defined by
\begin{align*}
\mathcal{B}_{\tau}[y]^{*}:&H^{1}(\Omega){\to} \left(L^{\infty}(\Omega)\right)',\\
{}&\langle \mathcal{B}_{\tau}^{*}[y]p,u\rangle_{(L^{\infty}(\Omega))', L^{\infty}(\Omega)}=\langle \mathcal{B}[u](y),p\rangle_{(H^{1}(\Omega))', H^{1}(\Omega)},
\end{align*}
which is also used in the derivation of the IDSM.

\subsection{The DSM in an operator form}\label{ssec:DSM}
We now reformulate the standard DSM for the abstract model \eqref{eqn1}. First we revisit the definitions of the background solution $y_\emptyset$ and the associated scattering field $y_*^s$ introduced in Section \ref{subsec_dsm} for the EIT problem.
For the general elliptic model \eqref{eqn1}, if the operator $\mathcal{A}$ is nonlinear, we define the background solution $y_{\emptyset}$ by
\begin{equation}\label{eqn24}
\mathcal{A}[y_{*}]y_{\emptyset} = f.
\end{equation}
Note that in \eqref{eqn24}, the operator $\mathcal{B}[u](y)$ does not appear and thus the inhomogeneity $u$ does not influence the background solution $y_\emptyset$ directly. However, since the function $y_{*}$ is unknown, in practice, $y_{*}$ has to be approximated by the potential corresponding to an estimate of $u$, denoted by $y(u)$.
Specifically, let $y(u)$ be the solution of \eqref{eqn1} with inhomogeneity $u$. Then the background solution $y_\emptyset$ can be approximated by $y_\emptyset(u)$, which solves
\begin{equation}\label{eqn:background-op}
\mathcal{A}[y(u)]y_\emptyset(u)=f,
\end{equation}
where the subscript $\emptyset$ indicates the background solution, including the dependence on
the inhomogeneity $u$ via $y(u)$.
This approximation leads to errors in the trace $\mathcal{T}y_\emptyset$  of the background solution $y_\emptyset$ and also the noisy scattering field $y^{s}_{d}$, since we typically do not have an accurate estimate $u$ in direct methods.
However, in the iterative scheme, if the estimate $u$ approaches the true inhomogeneity $u_{*}$, the linear model \eqref{eqn:background-op} will also more accurately approximate the model \eqref{eqn24}.
Like before, we define the exact scattering field $y^{s} (u) $ and noisy one $y^{s}_{d}(u)$ respectively by
\begin{equation}\label{eqn:scattering2}
y^{s}(u) = \mathcal{T}(y_\emptyset(u)-y_*)\quad \mbox{and} \quad y^{s}_{d} (u) = \mathcal{T}y_\emptyset(u)-y_{d},
\end{equation}
and moreover we use the shorthand notation $y_*^s = y^{s}(u_*)$ (which is consistent with that in Section \ref{subsec_dsm}). Note that for a linear model, the definition of $y_{\emptyset}(u)$ is independent of $u$, and thus the operator form of $y_{\emptyset}(u)$ covers also the linear model.

Next we recast the DSM in an abstract setting.
The construction of the index function $\eta$ of the DSM in \eqref{eqn2} involves two components: the numerator $\zeta$ and the denominator, denoted by $\mathcal{R}$.
The numerator $\zeta$ depends on the scattering field $y^{s}_*\equiv y^{s}(u_{*})$.
By combining the relation \eqref{eqn24} with the defining identity $\mathcal{A}[y_{*}]y_{*}+\mathcal{B}[u_{*}](y_{*})=f$, we deduce
\begin{equation*}
\mathcal{A}[y_{*}]\left(y_{*}-y_\emptyset(u_*)\right)=-\mathcal{B}[u_{*}](y_{*}):=-F.
\end{equation*}
Thus, the scattering field $y_*^s$ can be represented by
\begin{equation*}
y_*^s=\mathcal{T}\mathcal{A}[y_{*}]^{-1}F.
\end{equation*}
The forward operator $\mathcal{TA}[y_{*}]^{-1}$ maps the source $F$ to the scattering field $y_*^s$.
It can be represented in an integral form using Green's formula (see the identity \eqref{eqn:scattering}), while in an operator form it is given by $\mathcal{T}\mathcal{A}[y_{*}]^{-1}$.
Then $\zeta$ is computed by solving \eqref{eqn3}, in which $(-\Delta_{\Gamma})^\gamma y^{s}_{d}(u_{*})$ acts as a boundary source:
\begin{align*}
\langle \mathcal{A}[y_{*}]p,\zeta\rangle_{\left(H^{1}(\Omega)\right)', H^{1}(\Omega)}
=\left( \mathcal{T}p, (-\Delta_{\Gamma})^{\gamma}y^{s}_{d}(u_{*})\right)_{L^2(\Gamma)},\quad\forall p\in H^{1}(\Omega),
\end{align*}
or equivalently,
\begin{align*}
\langle \mathcal{A}[y_{*}]^{*}\zeta,p\rangle_{\left(H^{1}(\Omega)\right)', H^{1}(\Omega)}=\langle \mathcal{T}^{*}(-\Delta_{\Gamma})^{\gamma} y^{s}_{d}(u_{*}), p\rangle_{\left(H^{1}(\Omega)\right)', H^{1}(\Omega)},\quad\forall p\in H^{1}(\Omega).
\end{align*}
This is exactly the adjoint operator of the forward operator $\mathcal{T}\mathcal{A}[y_{*}]$.
In an operator form, it can be expressed by
\begin{align*}
\zeta =&\mathcal{A}[y_{*}]^{-*}\mathcal{T}^{*}(-\Delta_{\Gamma})^{\gamma}y^{s}_{d}(u_{*})\\
\approx&\mathcal{A}[y{_{*}}]^{-*}\mathcal{T}^{*}(-\Delta_{\Gamma})^{\gamma}y^{s}(u_{*})\\
=&\left(\mathcal{T}\mathcal{A}[y_{*}]^{-1}\right)^{*}(-\Delta_{\Gamma})^{\gamma}\left(\mathcal{T}\mathcal{A}[y_{*}]^{-1}\right)F.
\end{align*}
Note that for nonlinear models, the noisy scattering field $y^{s}_{d}(u_{*})$ and the background model $\mathcal{A}[y_{*}]$ are generally unavailable (since both $u_*$ and $y_*$ are unknown), and thus approximated by $y^{s}_{d}(u)$ and $\mathcal{A}[y(u)]$, respectively, where $u$ is an estimate of the unknown ground truth $u_{*}$.

The resolver $\mathcal{R}$ scales the index function $\eta$ pointwise, and acts as a multiplication operator, namely multiplying $\eta$ by the function $( G(x,\cdot), G(x,\cdot))_{H^{\gamma}(\Gamma)}^{-1}$ or its approximation.

In summary, the index function $\eta$ for the DSM is given by
\begin{align}
\eta =&\mathcal{R}\zeta=\mathcal{R}\left(\mathcal{T}\mathcal{A}[y(u)]^{-1}\right)^{*}(-\Delta_{\Gamma})^{\gamma}y^{s}_{d}(u)\nonumber\\
\approx&{\mathcal{R}\left(\mathcal{T}\mathcal{A}[y(u)]^{-1}\right)^{*}(-\Delta_{\Gamma})^{\gamma}y^{s}(u)\nonumber}\\
\approx&\mathcal{R}\left(\mathcal{T}\mathcal{A}[y_*]^{-1}\right)^{*}(-\Delta_{\Gamma})^{\gamma}y^{s}(u_*)\nonumber\\
=&\mathcal{R}\left(\mathcal{T}\mathcal{A}[y_{*}]^{-1}\right)^{*}(-\Delta_{\Gamma})^{\gamma}\left(\mathcal{T}\mathcal{A}[y_{*}]^{-1}\right)F.\label{eqn9}
\end{align}
In the derivation, the first approximation replaces the noisy data $y_{d}$ with the exact one $\mathcal{T}y_*$, whereas the second approximation is related to the nonlinearity of the operator $\mathcal{A}[y]$ and is actually an equality for linear models.
Here the source $F=\mathcal{B}_{\tau}[y_{*}]u_{*}$ encodes the inclusion $u_*$ to be identified.
It is connected with the scattering field $y^{s}(u_*)$ by the forward operator $\mathcal{T}\mathcal{A}[y_{*}]^{-1}$, i.e.,
$y^{s}(u_*)=\mathcal{T}\mathcal{A}[y_{*}]^{-1}F.$
Thus, the adjoint operator $(\mathcal{T}\mathcal{A}[y(u)]^{-1})^{*}$ can be used to lift $y^{s}_{d}(u_*)$ from the boundary $\Gamma$ into the domain $\Omega$.
To achieve the nearly mutual orthogonality, the preconditioner $(-\Delta_{\Gamma})^{\gamma}$ is employed, which leads to the kernel operator $\left(\mathcal{T}\mathcal{A}[y_*]^{-1}\right)^{*}(-\Delta_{\Gamma})^{\gamma}\left(\mathcal{T}\mathcal{A}[y_{*}]^{-1}\right)$.
An approximate inverse of the kernel operator is given by $\mathcal{R}$.

\subsection{The proposed IDSM}\label{ssec:IDSM}
In essence, the IDSM iterates the DSM process to refine the index function $\eta$, but with three key innovations to the DSM \eqref{eqn9}.
First, we incorporate the adjoint operator $\mathcal{B}_{\tau}$ into the forward operator $\mathcal{T}\mathcal{A}[y]^{-1}$. The operator $\mathcal{T}\mathcal{A}[y]^{-1}$ targets the source $F$ instead of the inclusion $u$.
By defining $\mathcal{G}[u]=\mathcal{T}\mathcal{A}[y(u)]^{-1}\mathcal{B}_{\tau}[y(u)]$, which maps from the set $S$ of admissible inclusions to the boundary data in $L^{2}(\Gamma)$, we obtain the forward operator for $u$:
\begin{equation}\label{eqn15}
y^{s} (u_{*})=\mathcal{T}\mathcal{A}[y_{*}]^{-1}F\quad \Rightarrow \quad y^{s}(u_{*})=\mathcal{G}[u_{*}]u_{*}.
\end{equation}
The operator $\mathcal{G}[u]$ directly relates the scattering field $y^{s}(u_*)$ (cf. \eqref{eqn:scattering2}) to the inhomogeneity $u_*$, and replaces $\mathcal{T}\mathcal{A}[y_{*}]^{-1}$ in \eqref{eqn9} in the IDSM.
Additionally, we replace the preconditioner $(-\Delta_{\Gamma})^{\gamma}$ with the regularized DtN map $\Lambda_{\alpha}(\mathcal{A}[y(u)])$ to improve numerical stability. The resolver $\mathcal{R}$ is initialized similarly to the DSM, using either $G(x,x')$ or its approximation, but modified iteratively using the low-rank correction technique. Below we elaborate these three steps.

First, given an initial guess $u_{0} = 0$, the index function $\eta_0$ in the IDRM is given by \eqref{eqn9}, i.e.,
\begin{equation}\label{eqn10}
\eta_{0}=\mathcal{R}_{0}\mathcal{G}[u_{0}]^{*}\Lambda_{\alpha}(\mathcal{A}[y(u_{0})])^{*}\Lambda_{\alpha}(\mathcal{A}[y(u_{0})])y^{s}_{d},
\end{equation}
where $\mathcal{R}_{0}$ is a multiplication operator that scales pointwise by the function $( G(x,\cdot), G(x,\cdot))_{H^{\gamma}(\Gamma)}^{-1}$ or its approximation.
Once $\eta_{0}$ is computed, it provides an initial estimate for $u$, which can be combined with additional prior information, e.g., box constraints, denoted by $u_{1} = \mathcal{P}(\eta_{0})$.

Then the IDSM proceeds to refine $\eta$ iteratively. Instead of simply replacing $u_{0}$ in \eqref{eqn10} by $u_{1}$, we also update the resolver $\mathcal{R}_{0}$.
In the original DSM, $\mathcal{R}$ is a pointwise multiplication operator, which is computationally efficient but has limited accuracy.
Also Green's function $G(x,x')$ may be not computed explicitly.
Thus, it is important to construct the resolver $\mathcal{R}$ numerically, which must approximate the inverse of the operator $G[u^*]$.
To motivate the construction, we assume that there is an auxiliary scattering field $\hat{y}^{s}$ associated with an estimate $\hat{u}$ such that
\begin{equation*}
\hat{y}^{s} = \mathcal{T}\mathcal{A}[y(\hat{u})]^{-1}\mathcal{B}_{\tau}[y(\hat{u})]\hat{u}.
\end{equation*}
Then following the preceding argument, the numerator  $\hat{\zeta}$ can be computed from $\hat y^s$ by
\begin{equation*}
\hat{\zeta} = \mathcal{B}_{\tau}[y(\hat{u})]^{*}\mathcal{A}[y(\hat{u})]^{-*}\hat{y}^{s}.
\end{equation*}
That is, there holds the identity $\hat{\zeta} = \mathcal{G}[\hat{u}]\hat{u}$.
Given an approximate inverse $\widehat{\mathcal{R}}$ of $\mathcal{G}[u_{*}]$ and
an estimate $\hat{u}$ of $u_{*}$, the auxiliary pair $(\hat{u}, \hat{\zeta})$ allows refining $\widehat{\mathcal{R}}$ by enforcing the relation $\widehat{\mathcal{R}}\hat{\zeta} = \hat{u}$.
These discussions motivate adjusting $\mathcal{R}_{0}$ using the auxiliary scattering field $\mathcal{T}(y_{\emptyset}(u_{1})-y(u_{1}))$ corresponding to $u_{1}$. Indeed, let $y(u_{1})$ be the solution of \eqref{eqn1} with the inclusion $u_{1}$:
\begin{equation*}
\mathcal{A}{[y(u_{1})]}y(u_{1})+\mathcal{B}[u_{1}](y(u_{1}))=f.
\end{equation*}
The relationship between $\mathcal{T}(y_{\emptyset}(u_{1})-y(u_{1}))$ and $u_{1}$ can be expressed using $\mathcal{G}[u]$ as
\begin{equation*}
\mathcal{T}(y_{\emptyset}(u_{1})-y(u_{1}))=\mathcal{T}\mathcal{A}[y(u_{1})]^{-1}\mathcal{B}_{\tau}[y(u_{1})]u_{1}{=\mathcal{G}[u_{1}]u_{1}}.
\end{equation*}
This formula follows exactly as for $u_{*}$ in \eqref{eqn15}.
Let
\begin{equation*}
\tilde{\zeta}_{1}=\mathcal{G}[u_{1}]\Lambda_{\alpha}(\mathcal{A}[y(u_{1})])^{*}\Lambda_{\alpha}(\mathcal{A}[y(u_{1})])\mathcal{T}\left[y_{\emptyset}(u_{1})-y(u_{1})\right].
\end{equation*}
Clearly, if the estimate $u_{1}$ were identical with the exact inclusion $u_*$, then $\tilde{\zeta}_{1}$ would serve as the ideal numerator of the index function $\eta$.
Thus  to obtain a more accurate resolver $\mathcal{R}_{1}$ than $\mathcal{R}_0$, we request that $\mathcal{R}_1$ satisfies the following relation:
\begin{align}\label{eqn:tangent} \mathcal{R}_{1}\tilde{\zeta}=u_{1}.
\end{align}
This choice ensures that the resolver $\mathcal{R}_{1}$ provides the exact prediction $u_{1}=\mathcal{R}_{1}\tilde{\zeta}_{1}$ for the problem with the ground truth $u_{1}$, which is expected to yield a more accurate indicator $\eta = \mathcal{R}_{1}\zeta_{1}$ for the problem. Numerically the relation \eqref{eqn:tangent} alone is insufficient to uniquely to determine $\mathcal{R}_1$. We draw inspiration from the optimization community: The relation \eqref{eqn:tangent} is formally identical with the secant equation for quasi-Newton type methods in the optimization literature. In the optimization community, researchers have proposed various ways to enforce the sequent equation, most prominently Broyden type low-rank corrections, including DFP (Davidon-Fletcher-Powell) update or BFG (Broyden-Fletcher-Goldfarb) update (see, e.g., \cite[p. 149]{Wright2006} and \cite[p. 169]{Schnabel1996}) for details). Therefore, we propose using these low-rank corrections to update the resolver $\mathcal{R}$.
Specifically, for any $\xi\in \left(L^{\infty}(\Omega)\right)'$, the updated part $\delta\mathcal{R}$ of the DFP and BFG corrections is implemented respectively as
\begin{align}
\delta\mathcal{R}_{\rm DFP}\xi& = \frac{\langle \xi,u_{k+1}\rangle}{\langle \tilde{\zeta}_{k+1},u_{k+1}\rangle}u_{k+1}
-\frac{\langle \xi,\mathcal{R}_{k}\tilde{\zeta}_{k+1}\rangle}{\langle \tilde{\zeta}_{k+1},\mathcal{R}_{k}\tilde{\zeta}_{k+1}\rangle}\mathcal{R}_{k}\tilde{\zeta}_{k+1},\label{eqn:DFP}\\
\delta\mathcal{R}_{\rm BFG}\xi& = \left(1+\frac{\langle \tilde{\zeta}_{k+1},\mathcal{R}_{k}\tilde{\zeta}_{k+1}\rangle}{\langle \tilde{\zeta}_{k+1},u_{k+1}\rangle}\right)\frac{\langle \xi,u_{k+1}\rangle}{\langle \tilde{\zeta}_{k+1},u_{k+1}\rangle}u_{k+1}\label{eqn:BFG}\\
&\qquad -\frac{\langle \xi,u_{k+1}\rangle}{\langle \tilde{\zeta}_{k+1},u_{k+1}\rangle}\mathcal{R}_{k}\tilde{\zeta}_{k+1}
-\frac{\langle \xi,\mathcal{R}_{k}\tilde{\zeta}_{k+1}\rangle}{\langle \tilde{\zeta}_{k+1},u_{k+1}\rangle}u_{k+1},\nonumber
\end{align}
where in the formulas, the notation $\langle\cdot, \cdot\rangle$ denotes the duality product $\langle\cdot, \cdot\rangle_{(L^{\infty}(\Omega))', L^{\infty}(\Omega)}$ between the space $L^\infty(\Omega)$ and  its dual $(L^\infty(\Omega))'$.
In practice, these corrections require computing only several duality products between $\xi$ and $u_{k+1}$ or $\mathcal{R}_{k}\tilde{\zeta}_{k+1}$, which is computationally efficient \cite[Chapter 6]{Wright2006} and also incurs a low-storage requirement.
Furthermore, either update maintains the symmetry and positive definiteness of the resolver $\mathcal{R}_{k}$.
This property is highly desirable since the resolver $\mathcal{R}_{k}$ approximates the inverse of a symmetric positive definite kernel operator $\left\{\Lambda_{\alpha}(\mathcal{A}[y_{*}])\mathcal{G}[u_{*}]\right\}^{*}\left\{\Lambda_{\alpha}(\mathcal{A}[y_{*}])\mathcal{G}[u_{*}]\right\}$.
In the context of optimization, \cite[Theorem 8.2.2]{Schnabel1996} states that for any objective function, regardless of whether they are quadratic or not, these two Broyden updates yield identical iterative sequences.
Thus, it is expected that the choice of the low-rank update rule should not affect much the computational results in the IDSM, and this is confirmed by the numerical experiments in Section \ref{sec_num}.

\begin{algorithm}[hbt!]
\caption{IDSM for nonlinear problems.}\label{alg1}
\begin{algorithmic}[1]
\State Initialize $\mathcal{R}_{0}$ and $u_{0}$.
\State Solve for $y_{0}=y(u_{0})$ from problem \eqref{eqn1} and set $y_{\emptyset}(u_{0})=\mathcal{T}\mathcal{A}[y_{0}]^{-1}f$.
\While{$k=0,1,\cdots,K-1$}
\State Solve $\zeta_{k}=\mathcal{G}[u_{k}]^{*}\Lambda_{\alpha}(\mathcal{A}[y_{k}])^{*}\Lambda_{\alpha}(\mathcal{A}[y_{k}])y^{s}_{d}(u_{k})$.
\State Compute $\eta_{k} = \mathcal{R}_{k}\zeta_{k}$.
\State Compute $u_{k+1}=\mathcal{P}(\eta_{k})$.
\State Solve for $y_{k+1}\equiv y(u_{k+1})$ from ~\eqref{eqn1} and $y_{\emptyset}(u_{k+1})=\mathcal{A}[y_{k+1}]^{-1}f$.
\State Solve $\tilde{\zeta}_{k+1}=\mathcal{G}[u_{k+1}]^{*}\Lambda_{\alpha}(\mathcal{A}[y_{k+1}])^{*}\Lambda_{\alpha}(\mathcal{A}[y_{k+1}])\mathcal{T}\left[y_{\emptyset}(u_{k+1})-y_{k+1}\right]$.
\State Calculate the update $\delta\mathcal{R}$ using either \eqref{eqn:DFP} or \eqref{eqn:BFG}.
\State Update $\mathcal{R}_{k+1}=\mathcal{R}_{k}+ \delta\mathcal{R}$.
\EndWhile
\end{algorithmic}
\end{algorithm}

We present the pseudocode for the IDSM in Algorithm \ref{alg1}.
In step 6, the operator $\mathcal{P}$ enforces additional prior information about the inhomogeneity $u_{k+1}\in S$, e.g., box constraints.
For the EIT problem in Example \ref{exam:EIT}, we may take the projection into the pointwise box constraint:
$$u_{k+1}=\mathcal{P}(\eta_{k})=\max({\eta_{k}, -\sigma_{0}+\varepsilon}),$$
where $\varepsilon>0$ is a small number, in order to maintain the positivity of $\sigma_{k+1}=u_{k+1}+\sigma_{0}$.
Here $\mathcal{R}_{k}$ consists of a multiplication operator complemented by about a few rank-one corrections, ensuring minimal storage requirement and efficient evaluation.
The primary workload lies in evaluating $\mathcal{G}[u_{k}]$, $\Lambda_{\alpha}(\mathcal{A}[y])$, along with their adjoints, necessitating solving several elliptic equations.
Despite the appearance of solving three elliptic equations, the computation of \eqref{eqn9} can be carried out by solving two elliptic problems only with the help of Lemma \ref{lemma2}:
First solve \eqref{eqn7} with $v=y^{s}_{d}$ and then solve \eqref{eqn11} with $v=z$. Then $\mathcal{B}_{\tau}[y(u)]w$ is the sought outcome.
\begin{lemma}\label{lemma2}
Let $(w,z)\in H^{1}(\Omega)\times L^{2}(\Gamma)$ be the solution of the following system:
\begin{equation}\label{eqn11}
\left\{
\begin{aligned}
{}\langle \mathcal{A}[y(u)]^{*}w,p\rangle_{(H^{1}(\Omega))', {H^{1}(\Omega)}}+( z,\mathcal{T}p)_{L^{2}(\Gamma)}&=0,&\quad\forall p\in H^{1}(\Omega),\\
-( \mathcal{T}w,q)_{L^{2}(\Gamma)}+( \alpha z,q)_{L^{2}(\Gamma)}&=(q,v)_{L^{2}(\Gamma)},&\quad\forall q\in L^{2}(\Gamma).
\end{aligned}\right.
\end{equation}
Then the following identity holds
$$\mathcal{G}[u]^{*}\Lambda_{\alpha}(\mathcal{A}[u])^{*}v=-\mathcal{B}_{\tau}[y(u)]^{*}w.$$
\end{lemma}
\begin{proof}
By the definition of the operator $\Lambda_{\alpha}(\mathcal{A}[y])v$ and \eqref{eqn11}, we know $z=\Lambda_{\alpha}(\mathcal{A}[y])^{*}v$.
Meanwhile \eqref{eqn11} states
\begin{equation*}
\langle \mathcal{A}[y(u)]w, p\rangle_{(H^{1}(\Omega))', {H^{1}(\Omega)}}+( z,\mathcal{T}p)_{L^{2}(\Gamma)}=0,\quad\forall p\in H^{1}(\Omega).
\end{equation*}
That is,
$w=-\mathcal{A}[y(u)]^{-*}\mathcal{T}^{*}z$, which directly leads to
\begin{equation*}
\mathcal{B}_{\tau}[y(u)]^{*}w=-\mathcal{B}_{\tau}[y(u)]^{*}\mathcal{A}[y(u)]^{-*}\mathcal{T}^{*}z=-\mathcal{G}[u]^{*}z=-\mathcal{G}[u]^{*}\Lambda_{\alpha}(\mathcal{A}[y(u)])^{*}v.
\end{equation*}
This completes the proof of the lemma.
\end{proof}

To enhance the efficiency of the IDSM for multiple pairs of Cauchy data, instead of applying the IDSM to each individual pair, one can aggregate the data using $\mathcal{B}_{\tau}[y]^{*}$.
Specifically, for each source $f_{\ell}$ and inhomogeneity $u$, let $y_\ell(u)$ be the solution to problem \eqref{eqn1}, and $y^{s}_{d,\ell}(u):=y_{\emptyset,\ell}(u)-y_{d,\ell}$ the corresponding noisy scattering field. The forward operator $\mathcal{T}\mathcal{A}[y_{\ell}(u)]^{-1}\mathcal{B}_{\tau}[y_{\ell}(u)]$ varies with the source $f_{\ell}$, and it is denoted by $\mathcal{G}_{\ell}[u]$. Nonetheless, the scattering field $y^{s}_{,\ell}(u):=y_{\emptyset,\ell}(u)-{\mathcal{T}y_{*,\ell}}$ results from the action of $\mathcal{G}_\ell[u]$ on the same $u$.
This suggests defining a composite forward operator $\mathcal{G}[u]:L^{\infty}(\Omega)\to L^{2}(\Gamma)^{L}$,
$v \,\mapsto \,\left(\mathcal{G}_{1}[u]v,\mathcal{G}_{2}[u]v,\cdots,\mathcal{G}_{L}[u]v\right)^\top$, by
\begin{equation*}
\left(y^{s}_{,1}(u_{*}),y^{s}_{,2}(u_{*}),\cdots,y^{s}_{,L}(u_{*})\right)^\top=\mathcal{G}[u_{*}]u_{*}.
\end{equation*}
Then the numerator $\zeta_{k}$ of the index function $\eta$ reads
\begin{equation}\label{eqn21}
\begin{aligned}
\zeta_{k}=&\mathcal{G}[u_{k}]^{*}\left\{\Lambda_{\alpha}(\mathcal{A}[y_{k,\ell}])^{*}\Lambda_{\alpha}(\mathcal{A}[y_{k,\ell}]) y_{d,\ell}^s(u)\right\}_{\ell = 1}^{L}\\
=&\sum_{\ell = 1}^{L}\mathcal{G}_{l}[u]^{*}\Lambda_{\alpha}(\mathcal{A}[y_{k,\ell}])^{*}\Lambda_{\alpha}(\mathcal{A}[y_{k,\ell}])y_{d,\ell}^s(u).
\end{aligned}
\end{equation}
Thus, the function $\zeta_{k}$ combines the information from all data pairs $\{(y_{d,\ell},f_{\ell})\}_{l=1}^{L}$.

Below we illustrate the IDSM on cardiac electrophysiology in Example \ref{exam:EC}.
\begin{example}[The IDSM for cardiac electrophysiology] First we discuss the case of one pair of Cauchy data. We initialize $u$ and $\mathcal{R}$ to $u_{0}=0$ and $\mathcal{R}_{0}=d(x,\Gamma)^{\gamma}$, respectively.
The first iteration $(k=0)$ begins with obtaining $y_{0}$ from \eqref{eqn:EC} with $\chi_{\omega}$ replaced by $u_{0}$ {\rm(}and $\chi_{\Omega\backslash\omega}$ by $1-u_{0}${\rm)}.
Then $y_{\emptyset}(u_{0})=\mathcal{T}y$ is computed by solving the following model $(k=0)$:
\begin{equation}\label{eqn25}
\left\{\begin{aligned}
-\Delta y+y_{k}^{2}y=&f,\quad\text{in }\Omega,\\
\frac{\partial}{\partial n}y=&0,\quad\text{on }\Gamma.
\end{aligned}\right.
\end{equation}
This yields $y^{s}_{d}(u_0):=y_{\emptyset}(u_{0})-y_{d}$, which is the first estimate of the scattering field $y^{s}(u^*)$.
The first iteration proceeds by computing $w$ as $$w=\mathcal{A}[y_{0}]^{-*}\mathcal{T}^{*}\Lambda_{\alpha}(\mathcal{A}[y(u_{0})])^{*}\Lambda_{\alpha}(\mathcal{A}[y(u_{0})])y^{s}_{d}(u_{0}).$$
In view of Lemma \ref{lemma2}, $w$ solves the following system $(k=0)$:
\begin{equation}\label{eqn12}
\left\{\begin{aligned}
-\Delta z+y_{k}^{2}z&=0,\quad \text{in }\Omega,\\
z+\alpha\frac{\partial}{\partial n}z&=y^{s}_{d}(u_{k}),\quad \text{on }\Gamma,
\end{aligned}\right.
\quad\mbox{and}\quad
\left\{\begin{aligned}
-\Delta w+y_{k}^{2}w&=0,\quad\text{in }\Omega,\\
\quad w+\alpha\frac{\partial}{\partial n}w&=\frac{\partial}{\partial n}z,\quad \text{on }\Gamma.
\end{aligned}\right.
\end{equation}
Then there holds the relation
\begin{equation*}
\frac{\partial }{\partial n}z=\Lambda_{\alpha}(\mathcal{A}[y_{0}])y^{s}_{d}(u_{0}).
\end{equation*}
This formula differs slightly from the definition of the regularized DtN map $\Lambda_{\alpha}(\mathcal{A}[y_{0}])$ in \eqref{eqn7}, but it is equivalent in view of Green's identity.
The numerator $\zeta_{0}$ is then computed as $\zeta_{0}=\mathcal{B}_{\tau}[y_{0}]^{*}w$.
Note that the operator $\mathcal{B}_{\tau}[y_{0}]^{*}\in L\left(H^{1}(\Omega),\left(L^{\infty}(\Omega)\right)'\right)$ is defined for any $p\in L^{\infty}(\Omega)$ as
\begin{align*}
{}&\langle \mathcal{B}_{\tau}[y_{0}]^{*}w, p\rangle_{(L^{\infty}(\Omega))', L^{\infty}(\Omega)}=\langle \mathcal{B}_{\tau}[y_{0}]p,w\rangle_{(H^{1}(\Omega))', H^{1}(\Omega)}\\
=&\langle \mathcal{B}[p](y_{0}), w\rangle_{(H^{1}(\Omega))', H^{1}(\Omega)}
=\int_{\Omega}p\left[(\sigma-1)\nabla y_{0}\cdot \nabla w-y_{0}^{3}w\right]\,\mathrm{d}x,
\end{align*}
in view of \eqref{eqn17} and the second formula of \eqref{eqn16}.
Thus, the function $\zeta_{0}=\mathcal{B}_{\tau}[y_{0}]^{*}w$ is given by
$$\zeta_{0}=(\sigma-1)\nabla y_{0}\cdot \nabla w-y_{0}^{3}w,$$
and when using $d(x,\Gamma)^{\gamma}$ to approximate $\|G(x,\cdot)\|_{H^\gamma(\Gamma)}^2$, the first indicator function $\eta_0$ is $$\eta_{0}(x)=\left(\mathcal{R}_{0}\zeta_{0}\right)(x)=d(x,\Gamma)^{\gamma}\zeta_{0}(x).$$
By incorporating the projection operator $\mathcal{P}$, the estimate of $u$ is updated to $u_{1}=\mathcal{P}(\eta_{0})$. In the context of CE, since $u$ is a characteristic function, the projection operator $\mathcal{P}$ is set such that $0\leq\mathcal{P}(\eta)\leq 1$ for all $\eta$. The potential $y_1$ is then given by $y_{1}=y(u_{1})$.
The model in \eqref{eqn25} is then recomputed with $y_{0}$ replaced by $y_{1}$ to obtain $y_{\emptyset}(u_{1})$.
Since $y_{1}$ is a {\rm(}hopefully{\rm)} better approximation of $y_{*}$ than $y_{0}$, $y^{s}_{d}(u_1)$ also approximates the scattering field $y^{s}(u_{*})$ better.
Moreover, $y_{\emptyset}(u_{1})-\mathcal{T}y_{1}$ is the scattering field for the inhomogeneity $u_{1}$.
By resolving the two equations in \eqref{eqn12} with the Neumann data $y_{\emptyset}(u_{1})-\mathcal{T}y_{1}$ and $y_{0}$ replaced by $y_{1}$, the numerator $\tilde{\zeta}_{1}$ is given by  $$\tilde{\zeta}_{1}=(\sigma-1)\nabla y_{1}\cdot \nabla w-y_{1}^{3}w.$$

Using the data pair $(u_{1},\tilde{\zeta}_{1})$, the resolver $\mathcal{R}_0$ is updated to $\mathcal{R}_{1}=\mathcal{R}_{0}+\delta \mathcal{R}$ using either \eqref{eqn:DFP} or \eqref{eqn:BFG}.
Then one proceeds to the next iteration.
At the $k$-th iteration, the procedure is the same as the first iteration, with $0$ replaced by $k$.
The only difference lies in the action of the resolver $\mathcal{R}_{k}$.
For example, with the DFP correction, $\eta_{1}=\mathcal{R}_{1}\zeta_{1}$ is calculated as:
\begin{equation*}
\begin{aligned}
\eta_{1}(x)=&\left(\mathcal{R}_{1}\zeta_{1}\right)(x)=\left(\mathcal{R}_{0}\zeta_{1}\right)(x)+\left(\delta \mathcal{R}\zeta_{1}\right)(x)\\
=&d(x,\Gamma)^{\gamma}\zeta_{1}(x)+\frac{\langle \zeta_{1}(x'),u_{1}(x')\rangle}{\langle \tilde{\zeta}_{1}(x'),u_{1}(x')\rangle}u_{1}(x)
-\frac{\langle \zeta_{1}(x'),d(x',\Omega)^{\gamma}\tilde{\zeta}_{1}(x')\rangle}{\langle \tilde{\zeta}_{1}(x'),d(x',\Gamma)^{\gamma}\tilde{\zeta}_{1}(x')\rangle}d(x,\Gamma)^{\gamma}\tilde{\zeta}_{1}(x),
\end{aligned}
\end{equation*}
where $\langle\cdot, \cdot\rangle$ denotes the duality product $\langle\cdot, \cdot\rangle_{(L^{\infty}(\Omega))', L^{\infty}(\Omega)}$ between the space $L^\infty(\Omega_{x'})$ and  its dual $(L^\infty(\Omega_{x'}))'$, and $d(x,\Omega)^{\gamma}\tilde{\zeta}_{1}(x)$ and $u_{1}(x)$ are stored  functions for the low-rank correction.

For multiple pairs of Cauchy data, we indicate a function corresponding to the $\ell$th pair of Cauchy data by the subscript $\ell$. At the $k$th iteration of the IDSM, each $w_{\ell}$ is computed using \eqref{eqn12} with the noisy scattering field $y_{d,\ell}^s(u_{k})$. Then by \eqref{eqn21}, the aggregated $\zeta_k$ is given by
\begin{equation*}
\zeta_{k}=\sum_{\ell=1}^{L}\left[(\sigma-1)\nabla y_{k,\ell}\cdot \nabla w_{\ell}-y_{k,\ell}^{3}w_{\ell}\right].
\end{equation*}
This represents the only change for multiple pairs.
Combining $\zeta_{k}$, $\mathcal{R}_{k}$ and $\mathcal{P}$ yields an estimate $u_{k+1}$. Then one computes the potential $y_{k+1,\ell}$ by:
\begin{equation*}
\left\{\begin{aligned}
-\nabla\cdot\left[\left(1+u_{k+1}(\sigma-1)\right)\nabla y_{k+1,\ell}\right]+(1-u_{k+1})y_{k+1,\ell}^{3}&=f_{\ell},&\text{ in }\Omega,\\
\frac{\partial}{\partial n}y_{k+1,\ell}&=0,&\text{ on }\Gamma.
\end{aligned}\right.
\end{equation*}
The process is repeated by solving \eqref{eqn12} with $y_{d,\ell}$ replaced by $\mathcal{T}y_{k+1,\ell}$, which leads to $\tilde{\zeta}_{k+1}$. Then we update $\mathcal{R}$ so that $\mathcal{R}_{k+1}$ satisfies
$\mathcal{R}_{k+1}\tilde{\zeta}_{k+1}=u_{k+1}$, which finishes the current iteration.
\end{example}

When the operator $\mathcal{A}$ is linear, the solution $y_{\emptyset}(u)$ of the background model is independent of the inclusion $u$. This allows computing certain quantities prior to the iterative process, thereby significantly reducing the overall computational expense.
For example, in EIT (with $\sigma_0\equiv1$), once the background solution $y_{\emptyset}(u)$ is computed, the preconditioner $\Lambda_{\alpha}(\mathcal{A})$ and the adjoint $\mathcal{A}^{-*}$ can be applied directly since they do not depend on $u$.
Specifically, two equations are solved sequentially (since for linear models, the quantity $y_{\emptyset}(u)$ does not require recomputation):
\begin{equation}\label{eqn23}
\left\{
\begin{aligned}
-\Delta z&=0,\quad \text{in }\Omega,\\
z+\alpha\frac{\partial}{\partial n}z&=y^{s}_{d}(u),\quad \text{on }\Gamma,\\
\end{aligned}\right.
\quad\mbox{and}\quad \left\{
\begin{aligned}
-\Delta w&=0,\quad \text{in }\Omega,\\
\quad w+\alpha\frac{\partial}{\partial n}w&=\frac{\partial}{\partial n}z,\quad \text{on }\Gamma.
\end{aligned}\right.
\end{equation}
Then at each iteration, computing {$\zeta_{k}$} requires only one pointwise multiplication: $$\zeta_{k}=\mathcal{B}_{\tau}[y_{k}]^*w=\nabla y_{k}\cdot\nabla w.$$
The details of the algorithm are given in Algorithm \ref{alg2}.
Each iteration of the IDSM requires solving three linear elliptic problems:
At step 7, $y_{k}$ is computed for the new inhomogeneity $u_{k}$;
at step 8, each one of the operator $\Lambda_{\alpha}(\mathcal{A})$ and $\mathcal{G}[u_{k+1}]^{*}\Lambda_{\alpha}(\mathcal{A})^{*}$ requires solving one elliptic problem.
Thus, it is more expensive than  one step of the standard regularized least-squares  of an elliptic inverse problem, which requires solving the primal problem and adjoint problem only for each step (and occasionally also Riesz map, e.g., when using the $H^1(\Omega)$ penalty).
However, the IDSM converges in a few steps, as shown in the  numerical experiments in Section \ref{sec_num}, which is attributed to the good initial value of $R_{0}$ provided by the standard DSM and the preconditioner $\Lambda_{\alpha}(\mathcal{A})$. Furthermore, in~\eqref{eqn23}, the two problems have different boundary conditions but share the same elliptic operator (or stiffness matrix in FEM implementation), which can be exploited to reduce the computational expense (e.g., assembling the stiffness matrix).

\begin{algorithm}
\caption{IDSM for linear problems.}\label{alg2}
\begin{algorithmic}[1]
\State Initialize $\mathcal{R}_{0}$ and $u_{0}$.
\State Solve $y_{\emptyset}{(u_{0})}=\mathcal{A}^{-1}f$ and $w=\mathcal{A}^{-*}\Lambda_{\alpha}(\mathcal{A})^{*}\Lambda_{\alpha}(\mathcal{A})y^{s}_{d}(u_{0})$.
\While{$k=0,1,\cdots,K-1$}
\State Compute $\zeta_{k}=\mathcal{B}_{\tau}[y_{k}]^{*}w$.
\State Compute $\eta_{k} = \mathcal{R}_{k}\zeta_{k}$.
\State Compute $u_{k+1}=\mathcal{P}(\eta_{k})$.
\State Solve $y_{k+1}=\left\{\mathcal{A}+\mathcal{B}[u_{k+1}]\right\}^{-1}f$.
\State Solve $\tilde{\zeta}_{k+1}=\mathcal{G}[u_{k+1}]^{*}\Lambda_{\alpha}(\mathcal{A})^{*}\Lambda_{\alpha}(\mathcal{A})\mathcal{T}\left[w-y_{k+1}\right]$.
\State Calculate the update $\delta\mathcal{R}$ using either \eqref{eqn:DFP} or \eqref{eqn:BFG}.
\State Update $\mathcal{R}_{k+1}=\mathcal{R}_{k}+\delta\mathcal{R}$.
\EndWhile
\end{algorithmic}
\end{algorithm}

Finally, note that if an approximation is employed for $\mathcal{R}_{0}$, e.g., using $d(x, \Gamma)^{\gamma}$, the scale of $\mathcal{R}_{0}$ may be incorrect, since the constant $C$ in \eqref{eqn22} may have been omitted. In Algorithms \ref{alg1} and \ref{alg2}, $\mathcal{R}_{0}$ undergoes low-rank correction and it is important to scale $\mathcal{R}_0$ properly.
A convenient choice is to scale $\mathcal{R}_{0}$ using the factor $\|u_{1}\|_{L^{1}(\Omega)}/\|\mathcal{R}_{0}\tilde{\zeta}_{1}\|_{L^{1}(\Omega)}$ prior to the low-rank update.

\section{Numerical experiments and discussions}\label{sec_num}

In this section we present numerical results to illustrate distinct features of the IDSM. All the experiments are conducted in an elliptical domain $\Omega=\{(x_1,x_2): \frac{x_1^{2}}{1} + \frac{x_2^2}{0.64}< 1\}$.
The exact data $y_{*}(x)$ is generated by solving problem \eqref{eqn1} on a fine mesh.
To generate noisy data $y_{d}$, we add random noise into $\mathcal{T}y_{*}$ pointwise by
\begin{equation*}
y_{d}(x)=y_{*}(x)+\varepsilon \delta(x)\left|y_{\emptyset}(0)(x)-y_{*}(x)\right|,
\end{equation*}
where $\delta(x)$ is uniformly distributed in $(-1, 1)$ and $\varepsilon$ refers to the relative noise level.
This choice follows that of \cite[(6.1)]{Zou2015} and \cite[(6.1)]{Zou2014}.
The codes for reproducing the numerical results will be made available at the following github link \url{https://github.com/RaulWangfr/IDSM-elliptic.git}.

The experiments contain five examples. Example 1 considers EIT at two noise levels ($\varepsilon = 10\%$ and $\varepsilon = 30\%$), in order to show the resilience of the IDSM with respect to data noise. Example 2 shows the capability of the IDSM to discern inclusions of different types, i.e., simultaneously reconstructing conductivity and potential.
Example 3 studies the impact of the regularization parameter $\alpha$ in the regularized DtN map $\Lambda_\alpha(\mathcal{A})$ on the performance of the IDSM in the context of diffuse optical tomography.
Example 4 illustrates the applicability of the IDSM to semilinear models. Finally, Example 5 examines a semilinar elliptic model featuring two inclusions with distinct intensities. This example shows that while the IDSM is unable to precisely recover the intensities of the inclusions, it is capable of distinguishing between inclusions with different coefficients.
These examples show the versatility of the IDSM for addressing a broad range of elliptic inverse problems.

\subsection{Example 1: electrical impedance tomography}

This experiment considers the EIT problem.
The conductivity of the background medium is $1$, and the region enclosed by a white square, which represents the inclusion, has a conductivity $0.3$ {(i.e., $u=-0.7$)}.
To test the resilience of the IDSM against noise, we take $\varepsilon=30\%$ and $\varepsilon=10\%$. Two pairs of Cauchy data are utilized, corresponding to the external boundary source $f_{1}=x_{1}$ and $f_{2}=x_{2}$.
We apply Algorithm \ref{alg2} with $\alpha=1.0$ in the regularized DtN map $\Lambda_\alpha{(\mathcal{A})}$. The resolver
$\mathcal{R}_{0}$ is initialized to $\|\nabla\Phi_{x}\|_{L^{2}(\Gamma)}$ with
$\Phi_{x}(x^{\prime})=-\frac{1}{2\pi}\ln(|x-x^{\prime}|)$
being the fundamental solution.
For the operator $\mathcal{P}$, we employ the pointwise projection to {the} box constraint $
\mathcal{P}(\eta)=\max(\min(\eta, 1.0), 0.01).$

Fig. \ref{fig1} shows the estimate $u_{k}$ during the iteration, where the squares indicate the inclusion locations.
The initial estimate $u_1$ ($k=1$) is inaccurate due to the presence of noise.
By the $6$th and $11$th iterations, the estimates $u_6$ and $u_{11}$ have noticeably improved the accuracy for localizing the inclusions.
The recovered regions become more concentrated and distinct, indicating that the IDSM has effectively converged to the exact inclusion locations.
The IDSM remains stable for up to $\varepsilon=30\%$ noise, and the recovered results are largely comparable with that for $\varepsilon=10\%$, clearly showing the robustness of the IDSM. Moreover, the two correction schemes work equally well.

\newcommand{\figlen}{0.28\textwidth}
\begin{figure}[hbt!]
\centering
\begin{tabular}{ccc}
\includegraphics[width = \figlen, trim = {0.0cm 0.0cm 0.0cm 1.5cm}, clip]{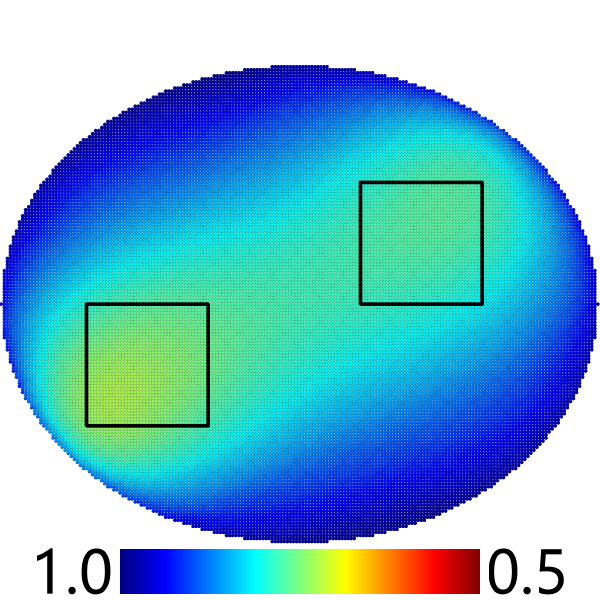}&
\includegraphics[width = \figlen, trim = {0.0cm 0.0cm 0.0cm 1.5cm}, clip]{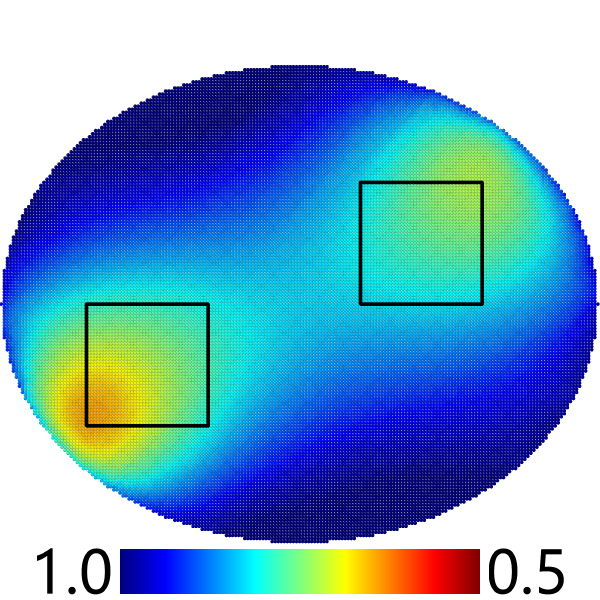}&
\includegraphics[width = \figlen, trim = {0.0cm 0.0cm 0.0cm 1.5cm}, clip]{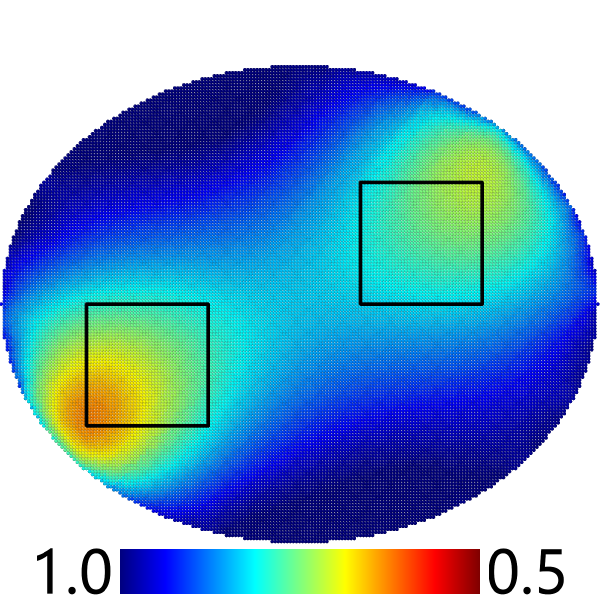}\\
\includegraphics[width = \figlen, trim = {0.0cm 0.0cm 0.0cm 1.5cm}, clip]{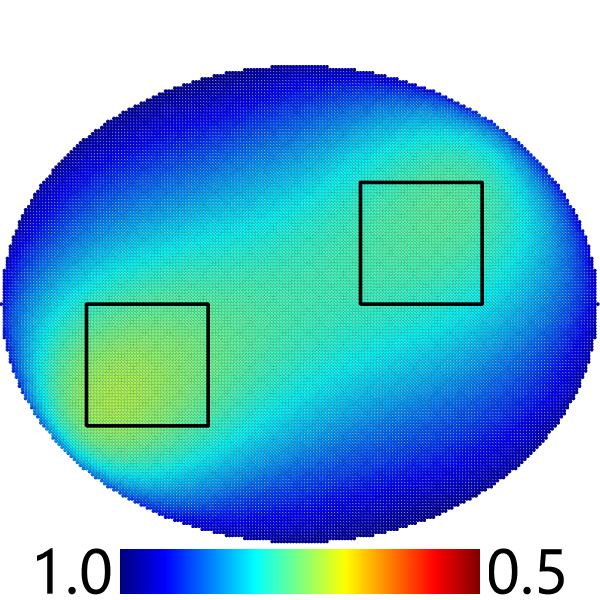}&
\includegraphics[width = \figlen, trim = {0.0cm 0.0cm 0.0cm 1.5cm}, clip]{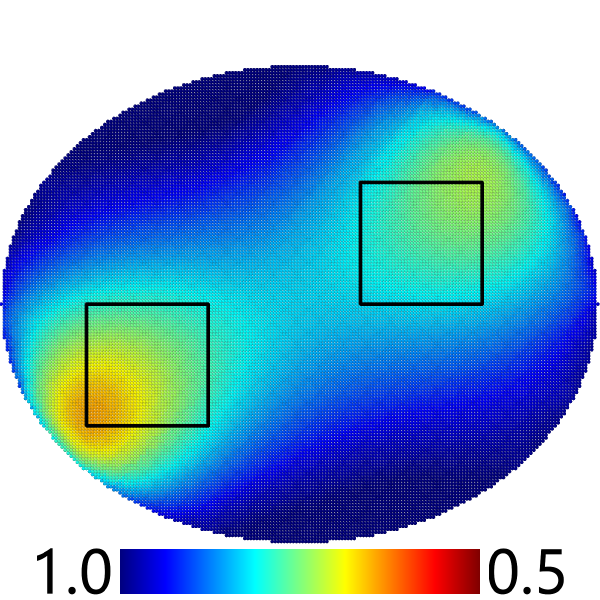}&
\includegraphics[width = \figlen, trim = {0.0cm 0.0cm 0.0cm 1.5cm}, clip]{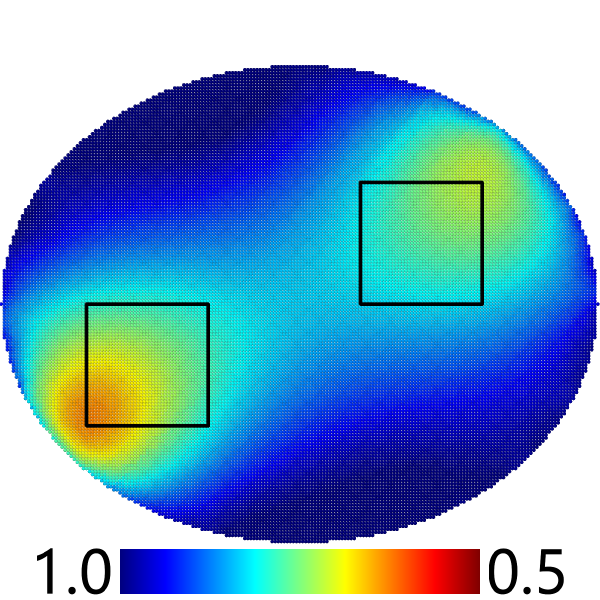}\\
\includegraphics[width = \figlen, trim = {0.0cm 0.0cm 0.0cm 1.5cm}, clip]{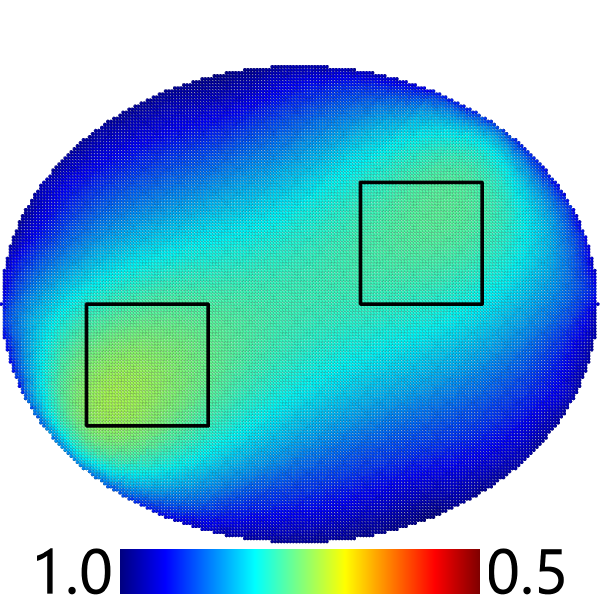}&
\includegraphics[width = \figlen, trim = {0.0cm 0.0cm 0.0cm 1.5cm}, clip]{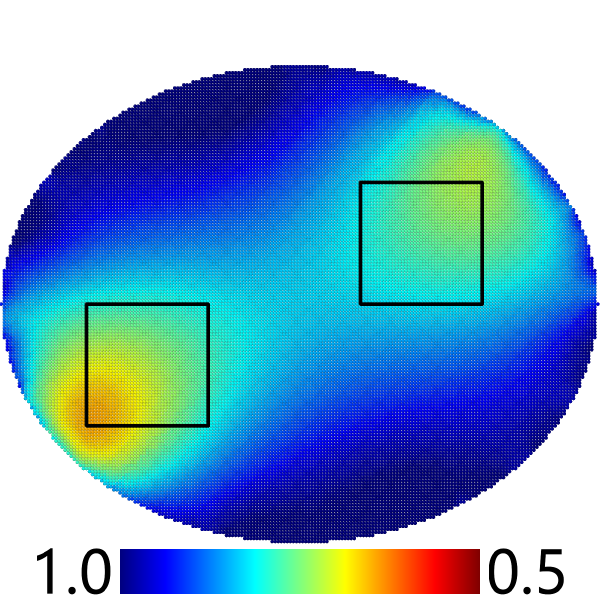}&
\includegraphics[width = \figlen, trim = {0.0cm 0.0cm 0.0cm 1.5cm}, clip]{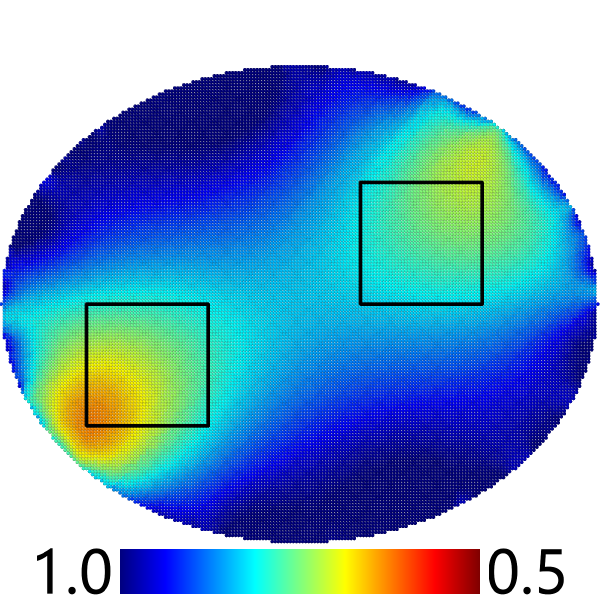}\\
\includegraphics[width = \figlen, trim = {0.0cm 0.0cm 0.0cm 1.5cm}, clip]{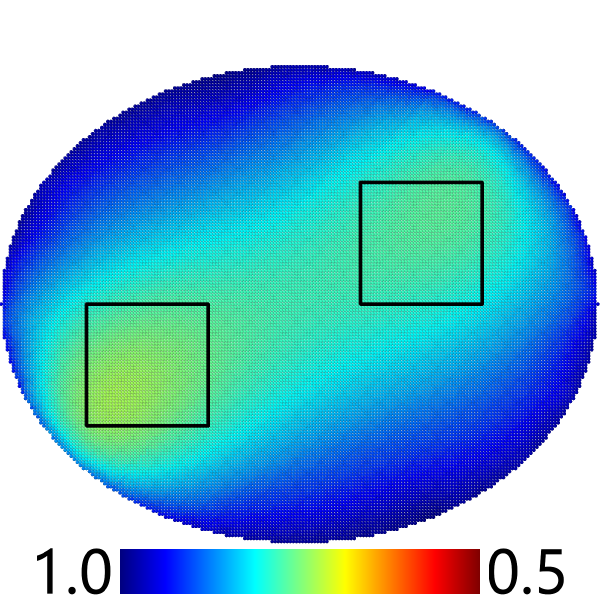}&
\includegraphics[width = \figlen, trim = {0.0cm 0.0cm 0.0cm 1.5cm}, clip]{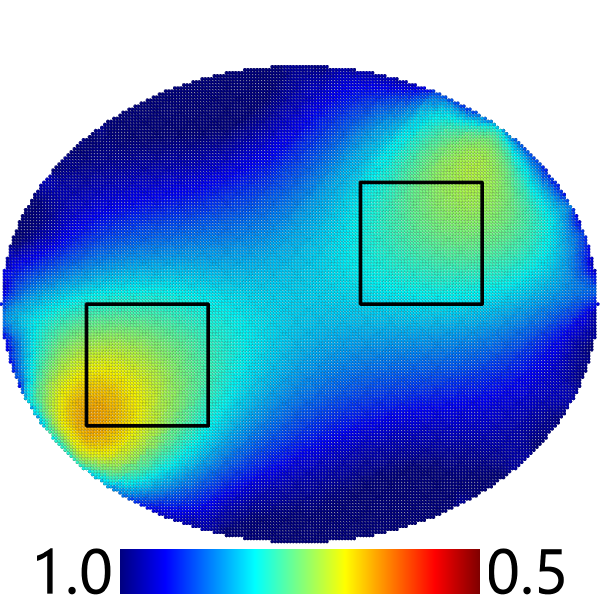}&
\includegraphics[width = \figlen, trim = {0.0cm 0.0cm 0.0cm 1.5cm}, clip]{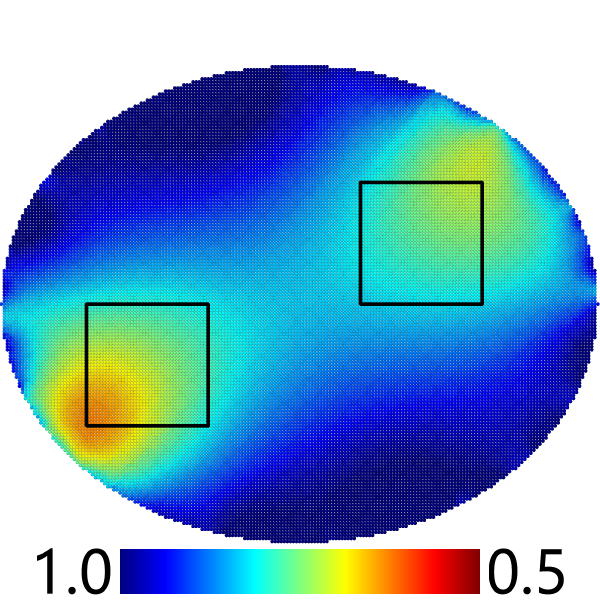}\\
(a) $k=1$ & (b) $k=6$ & (c) $k=11$
\end{tabular}
\caption{Visualization of $u_{k}$ for Example 1.
The conductivity is $0.3$ inside the inclusions and $1.0$ outside.
The first and second rows show the estimates for $\varepsilon = 10\%$, by the BFG and DFP corrections, respectively.
The third and fourth rows depict the estimates for $\varepsilon = 30\%$, by the BFG and DFP corrections, respectively.\label{fig1}}
\end{figure}

\subsection{Example 2: simultaneous recovery of conductivity and potential}

This example showcases the capability of the IDSM to differentiate inhomogeneities of distinct physical properties.
The governing equation of the model reads
\begin{equation*}
\left\{
\begin{aligned}
-\nabla \cdot \left[\left(1+u_{c}\right)\nabla y\right]+(1+u_{v})y&=0,&\text{ in }\Omega,\\
\frac{\partial}{\partial n}y&=f,&\text{ on }\Gamma,
\end{aligned}
\right.
\end{equation*}
and we aim to simultaneously recover the locations of conductivity inclusion $u_{c}$ and potential inclusion $u_{v}$.
The homogeneous operator $\mathcal{A}$ of the problem is given by
\begin{equation*}
\langle \mathcal{A}y, w\rangle_{(H^{1}(\Omega))', H^{1}(\Omega)}=\int_{\Omega}\nabla y\cdot\nabla w + yw\,{\rm d}x,
\end{equation*}
and the interaction of the inhomogeneities $[u_{c},u_{v}]$ with the governing equation is specified through
\begin{equation*}
\langle \mathcal{B}[u_{c}, u_{p}](y), w\rangle_{(H^{1}(\Omega))', H^{1}(\Omega)}=\int_{\Omega}u_{c}\nabla y\cdot\nabla w + u_{p}yw\,{\rm d}x.
\end{equation*}
In the figures, the region enclosed by the white square is the support of $u_{c}$, inside which the conductivity value is $0.3$ (i.e., $u_{c}=-0.7$), and the region enclosed by the black square is the support of $u_{v}$, inside which the potential value is $6$ (i.e., $u_{v}=5$).
In the experiment, $\varepsilon=10\%$ noise is added to the data, and two pairs of Cauchy data are utilized, corresponding to external boundary source $f_{1}=x_{1}$ and $f_{2}=x_{2}$.
The parameter $\alpha$ in the regularized DtN map $\Lambda_\alpha{(\mathcal{A})}$ is set to $0.1$.
Algorithm \ref{alg2} is applied with $\mathcal{R}_{0}=d(x,\Gamma)^{2}$. For $\mathcal{P}$, like before, we employ the projection to box constraint: $
\mathcal{P}([\eta_{c}, \eta_{p}])=[\max(\min(\eta_{c}, 1.0), 0.01), \max(\min(\eta_{p}, 10.0), 1.0)].$

\begin{figure}[hbt!]
\centering
\begin{tabular}{ccc}
\includegraphics[width = \figlen, trim = {0.0cm 0.0cm 0.0cm 1.5cm}, clip]{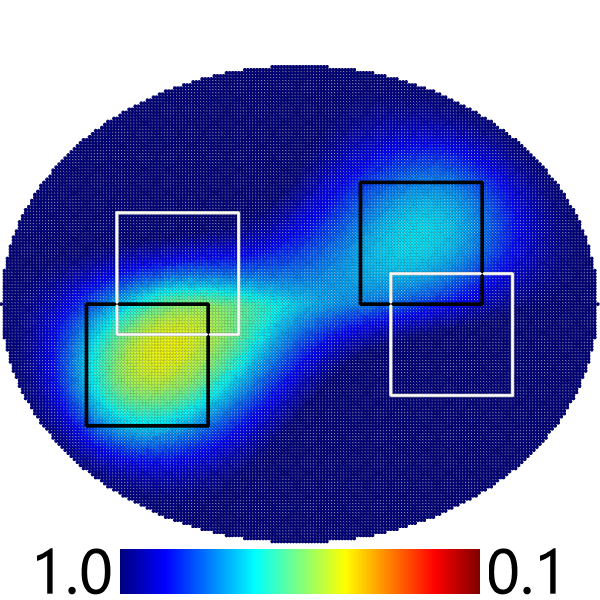}
&\includegraphics[width = \figlen, trim = {0.0cm 0.0cm 0.0cm 1.5cm}, clip]{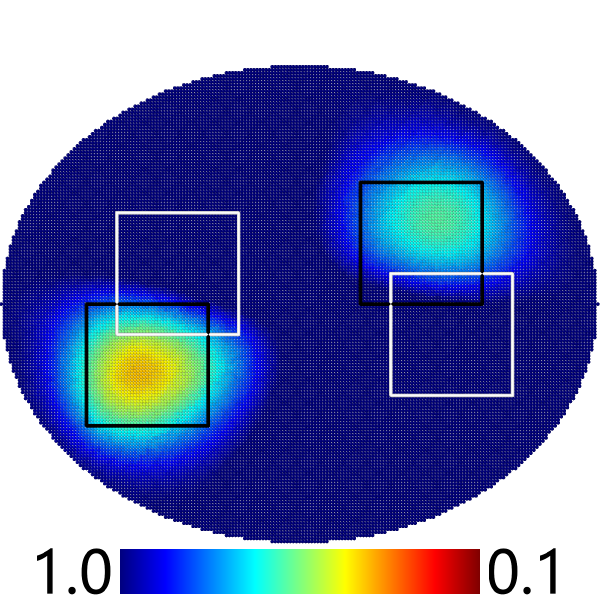}
&\includegraphics[width = \figlen, trim = {0.0cm 0.0cm 0.0cm 1.5cm}, clip]{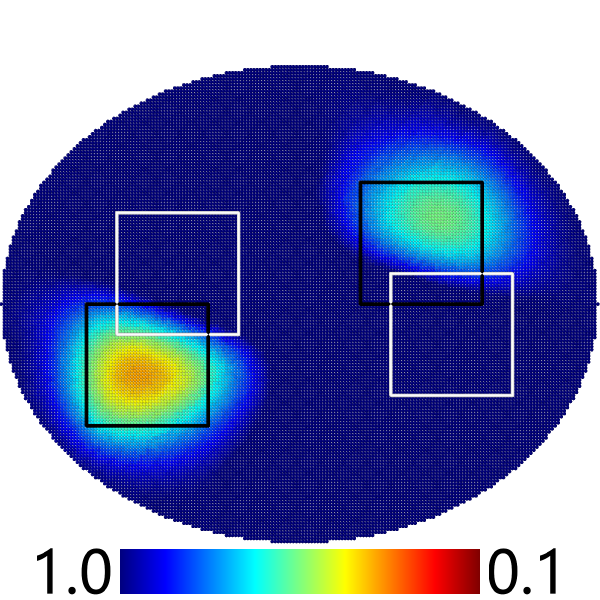}\\
\includegraphics[width = \figlen, trim = {0.0cm 0.0cm 0.0cm 1.5cm}, clip]{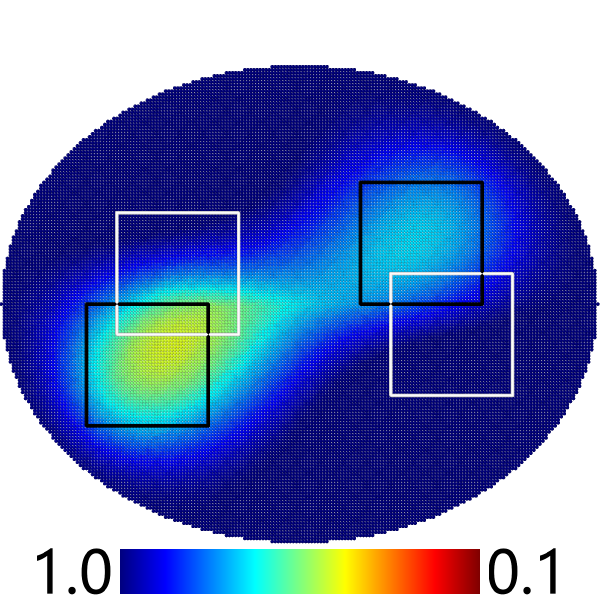}
&\includegraphics[width = \figlen, trim = {0.0cm 0.0cm 0.0cm 1.5cm}, clip]{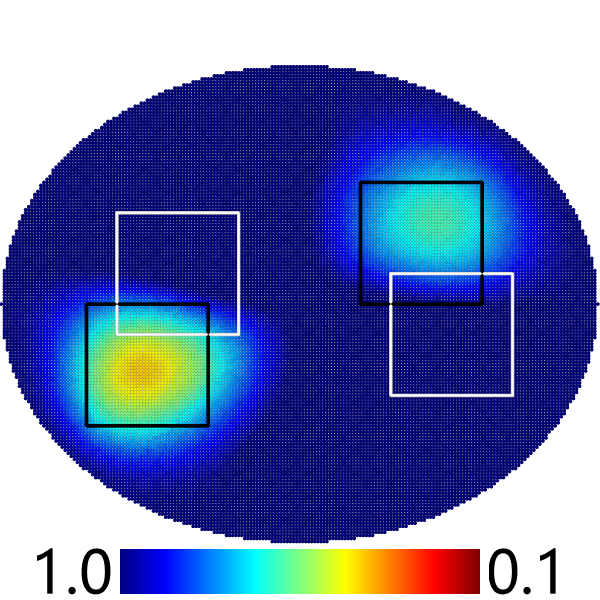}
&\includegraphics[width = \figlen, trim = {0.0cm 0.0cm 0.0cm 1.5cm}, clip]{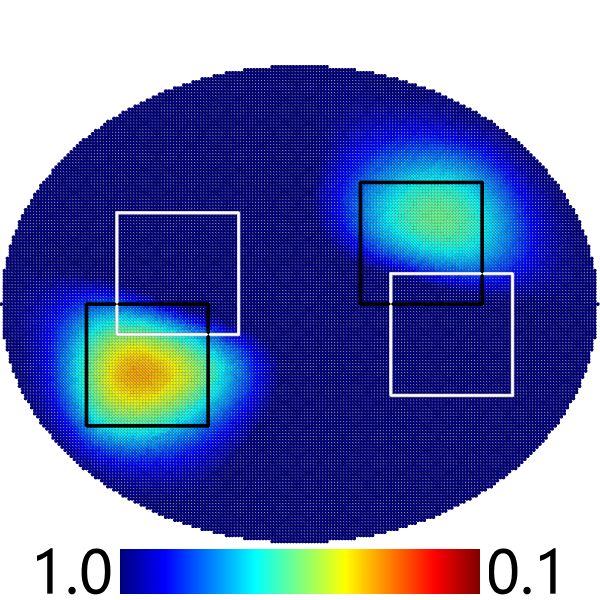}\\
\includegraphics[width = \figlen, trim = {0.0cm 0.0cm 0.0cm 1.5cm}, clip]{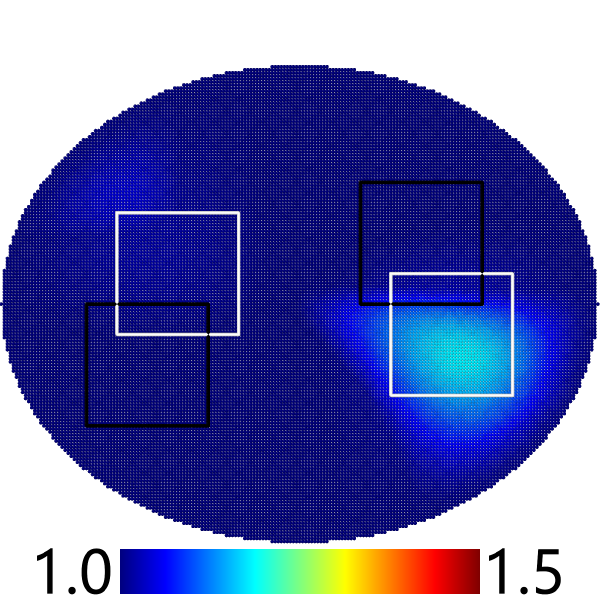}
&\includegraphics[width = \figlen, trim = {0.0cm 0.0cm 0.0cm 1.5cm}, clip]{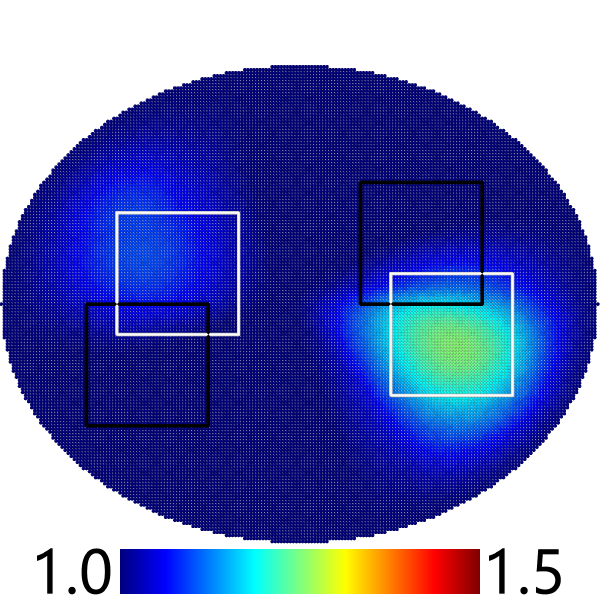}
&\includegraphics[width = \figlen, trim = {0.0cm 0.0cm 0.0cm 1.5cm}, clip]{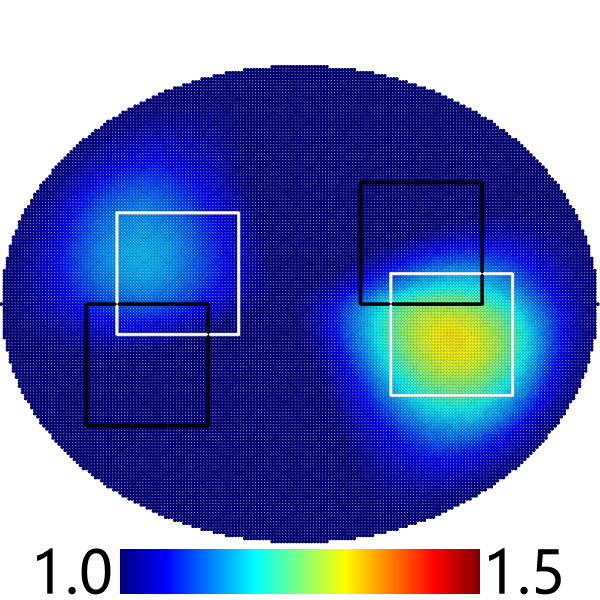}\\
\includegraphics[width = \figlen, trim = {0.0cm 0.0cm 0.0cm 1.5cm}, clip]{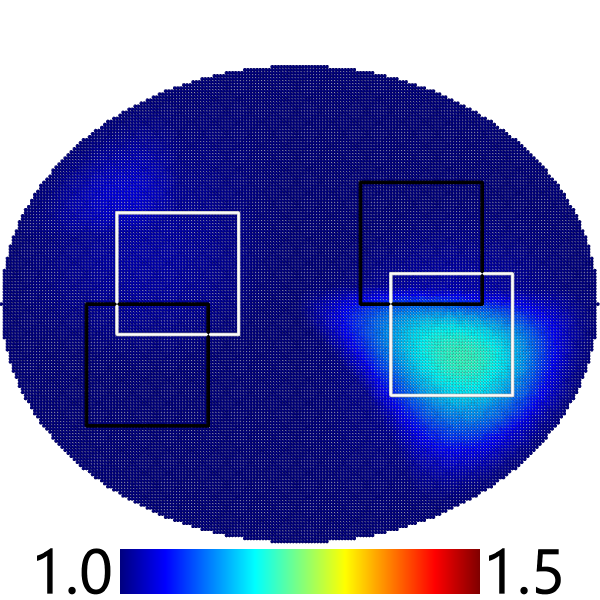}
&\includegraphics[width = \figlen, trim = {0.0cm 0.0cm 0.0cm 1.5cm}, clip]{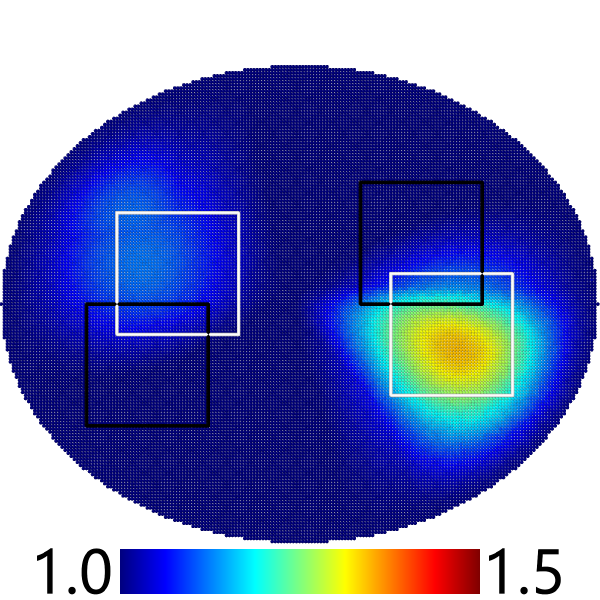}
&\includegraphics[width = \figlen, trim = {0.0cm 0.0cm 0.0cm 1.5cm}, clip]{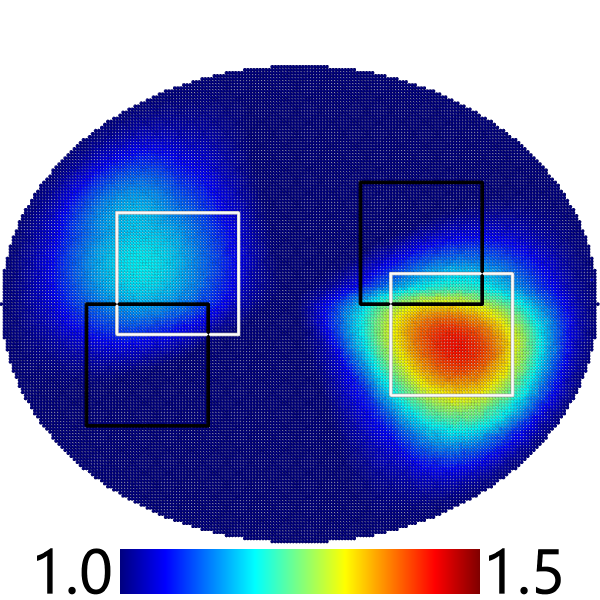}\\
(a) $k=1$ & (b) $k=6$ & (c) $k=11$
\end{tabular}
\caption{Visualization of $u_{k}$ for Example 2 with a noise level $\varepsilon = 10\%$. The conductivity is $0.3$ inside the conductivity inclusions and $1.0$ outside, and the potential is $6$ inside the potential inclusions and $1.0$ outside.
The first two rows show the estimates for conductivity inclusion, by the BFG and DFP corrections, respectively.
The last two rows depict the estimates for the potential inclusion, by the BFG and DFP corrections, respectively.
}
\label{fig2}
\end{figure}

The numerical results in Fig. \ref{fig2} show that the initial estimate $u_1$ ($k=1$) provides a baseline for the IDSM, which is then refined by successive iterations.
By the $6$th iteration, the estimate $u_6$ by the IDSM becomes more localized within the regions of interest and can more clearly distinguish the two types of inclusions. By the $11$th iteration, the method has effectively converged, and the reconstruction provides an accurate recovery of the inclusion locations.
This shows the robustness and versatility of the IDSM in handling multiple inhomogeneities of different physical properties.

\subsection{Example 3: diffuse optical tomography}

This example investigates the performance of the IDSM on diffuse optical tomography \cite{Arridge:1999}.
The problem is governed by the following model:
\begin{equation*}
\left\{
\begin{aligned}
-\Delta y+uy&=0,&\text{ in }\Omega,\\
\frac{\partial}{\partial n}y&=f,&\text{ on }\Gamma.
\end{aligned}
\right.
\end{equation*}
Thus the associated operators $\mathcal{A}$ and $\mathcal{B}$ are respectively defined by
\begin{align*}
\langle \mathcal{A}y, w\rangle_{\left(H^{1}(\Omega)\right)', H^{1}(\Omega)}&=\int_{\Omega}\nabla y\cdot\nabla w\,{\rm d}x,\\
\langle \mathcal{B}[u](y), w\rangle_{\left(H^{1}(\Omega)\right)', H^{1}(\Omega)}&=\int_{\Omega}uyw\,{\rm d}x.
\end{align*}
In the figures, the region enclosed by a black square indicates the support of the inclusion $u$, in which the potential value is 6.
The experiment is carried out for noisy data with a noise level $\varepsilon=10\%$.
We use one pair of Cauchy data, corresponding to the boundary source $f=x_{1}$.
Algorithm \ref{alg2} is applied with the resolver $\mathcal{R}$ initialized to $\|\nabla\Phi_{x}\|^{3}_{H^{1}(\Gamma)}$, where
$\Phi_{x}(x^{\prime})=-\frac1{2\pi}\ln(|x-x^{\prime}|)$
is the fundamental solution of $\mathcal{A}$ in $\mathbb{R}^{2}$.
The operator $\mathcal{P}$, like before, is the pointwise projection operator
$\mathcal{P}(\eta)=\max(\min(\eta, 10.0), 0.0)$.

\begin{figure}[hbt!]
\centering
\begin{tabular}{ccc}
\includegraphics[width =\figlen, trim = {0.0cm 0.0cm 0.0cm 1.5cm}, clip]{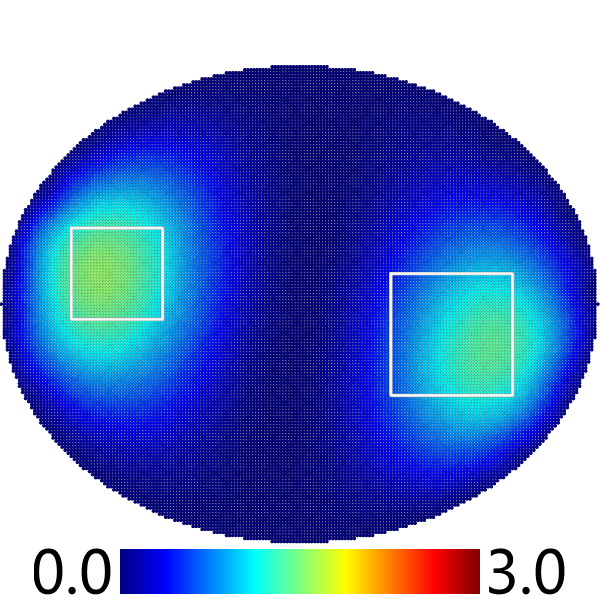}
&\includegraphics[width = \figlen, trim = {0.0cm 0.0cm 0.0cm 1.5cm}, clip]{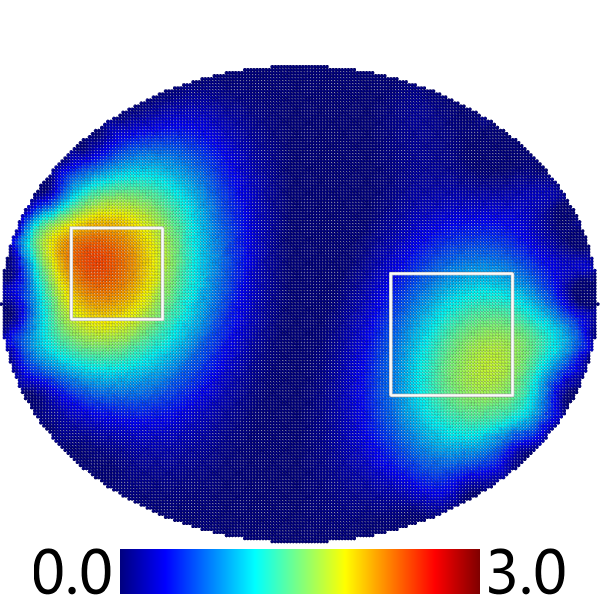}
&\includegraphics[width = \figlen, trim = {0.0cm 0.0cm 0.0cm 1.5cm}, clip]{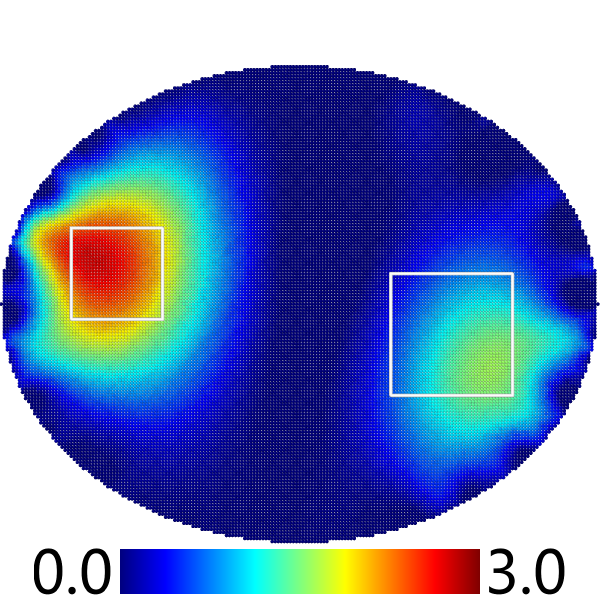}\\
\includegraphics[width = \figlen, trim = {0.0cm 0.0cm 0.0cm 1.5cm}, clip]{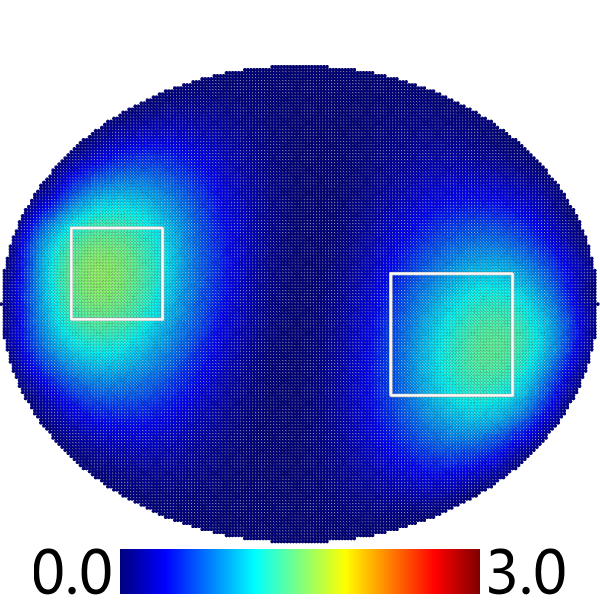}
&\includegraphics[width = \figlen, trim = {0.0cm 0.0cm 0.0cm 1.5cm}, clip]{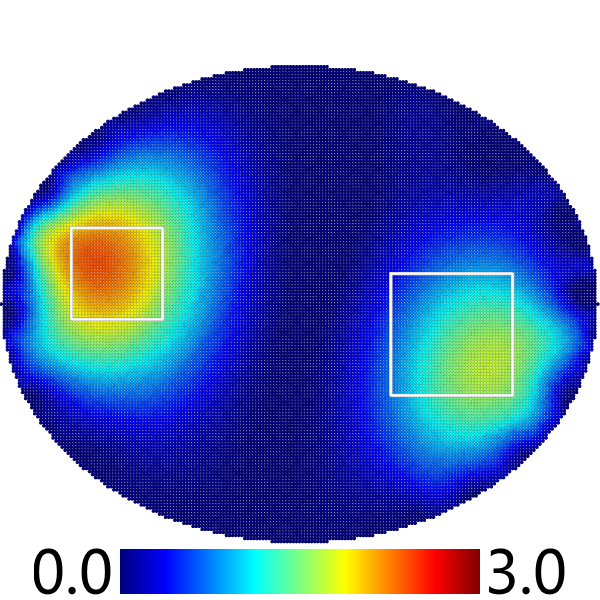}
&\includegraphics[width = \figlen, trim = {0.0cm 0.0cm 0.0cm 1.5cm}, clip]{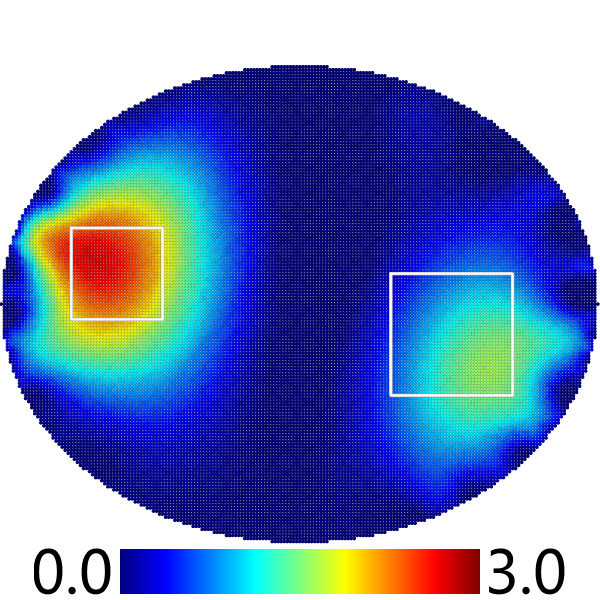}\\
(a) $k=1$ & (b) $k=11$ & (c) $k=21$\\
\includegraphics[width = \figlen, trim = {0.0cm 0.0cm 0.0cm 1.5cm}, clip]{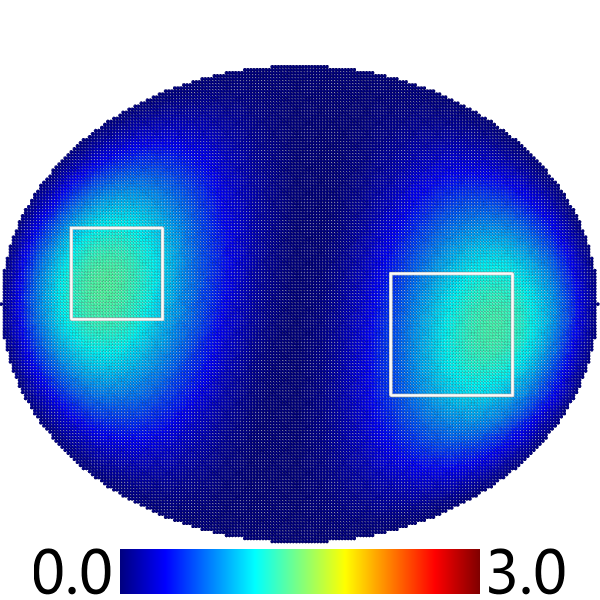}
&\includegraphics[width = \figlen, trim = {0.0cm 0.0cm 0.0cm 1.5cm}, clip]{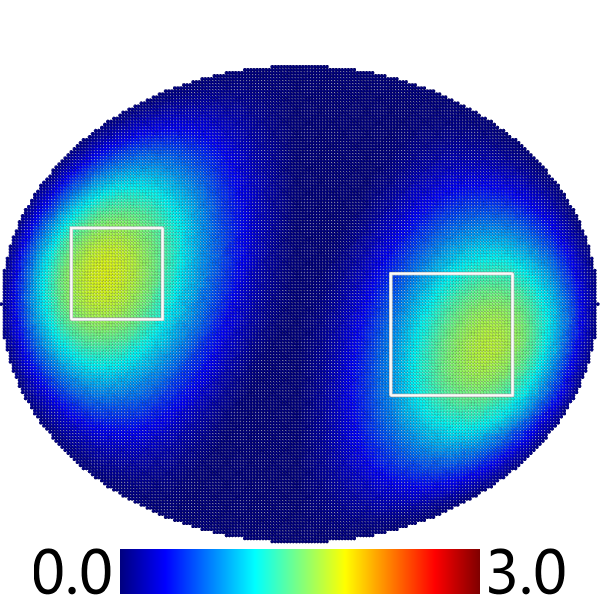}
&\includegraphics[width = \figlen, trim = {0.0cm 0.0cm 0.0cm 1.5cm}, clip]{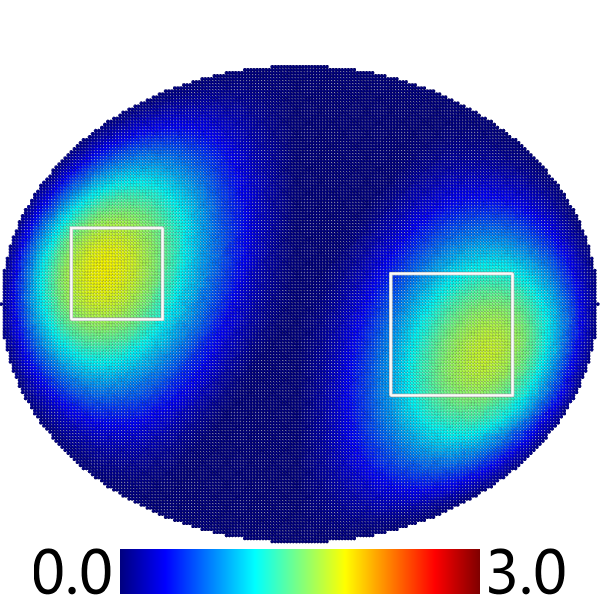}\\
\includegraphics[width = \figlen, trim = {0.0cm 0.0cm 0.0cm 1.5cm}, clip]{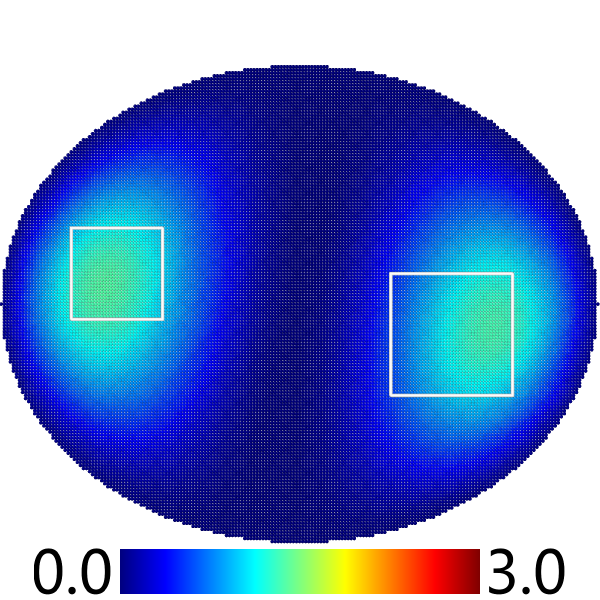}
&\includegraphics[width = \figlen, trim = {0.0cm 0.0cm 0.0cm 1.5cm}, clip]{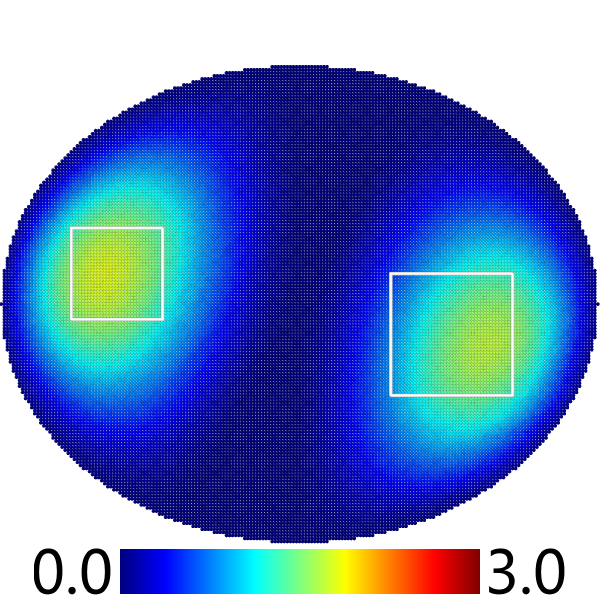}
&\includegraphics[width = \figlen, trim = {0.0cm 0.0cm 0.0cm 1.5cm}, clip]{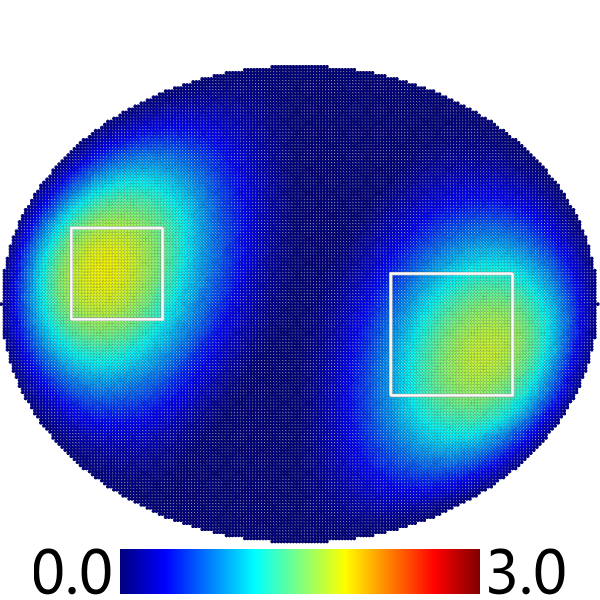}\\
{(d)} $k=1$ & {(e)} $k=6$ & {(f)} $k=11$
\end{tabular}
\caption{Visualization of $u_{k}$ for Example 3
with a noise level $\varepsilon = 10\%$. The potential is $6$ inside the  inclusions and $0.0$ outside.
The first and second rows show the estimates with $\alpha=\text{1e-3}$, using the BFG and DFP corrections, respectively.
The third and fourth rows depict the estimates with $\alpha=1$, using the BFG and DFP corrections, respectively.\label{fig3}}
\end{figure}

Fig. \ref{fig3} presents the reconstruction results for two choices of the penalty parameter $\alpha$ in the regularized DtN map $\Lambda_\alpha(\mathcal{A})$: $\alpha=\text{1e-3}$ and $\alpha=$1 .
The first two rows are for the choice $\alpha=\text{1e-3}$, which is an under-regularized scenario.
The estimates $u_k$ at $k=11$ and $k=21$ are only slightly affected by the presence of data noise.
Nonetheless, the IDSM maintains excellent stability over an extended number of iterations and delivers results that are sufficiently reliable.
The third and fourth rows are for $\alpha=1$ in $\Lambda_\alpha(\mathcal{A})$, which is an over-regularized scenario.
The initial estimate $u_1$ ($k=1$) cannot accurately capture the inclusions.
However, as the iteration proceeds to the $6$th and $11$th steps, the IDSM manages to refine the resolution of the estimate $u_k$, slowly but correctly resolving the exact inclusion locations.
These observations highlight the robustness of the IDSM with respect to the choice of the penalty parameter $\alpha$.
While the choice of $\alpha$ cannot fully filter out the noise or may impact the initial iterations, the IDSM enjoys high robustness, and its iterative nature allows continuous improvement, leading to an accurate localization of inclusions.
This feature is especially beneficial in real-world applications where determining an optimal penalty parameter $\alpha$ is challenging.

\subsection{Example 4: cardiac electrophysilogy}

This example involves a semilinear model in cardiac electrophysiology described in Example \ref{exam:EC}.
The conductivity $\sigma$ within the inhomogeneity (which physically represents the ischemia) is set to $\text{1e-4}$, consistent with that in \cite{beretta2018a}.
Algorithm \ref{alg1} is employed with $u$ initialized to $0$. The resolver $\mathcal{R}$ is initialized to the distance function $d(x,\Gamma)$, and the map $\mathcal{P}$ is given by
$\mathcal{P}(\eta_{k+1})=0.8u_{k}+0.2\frac{\eta_{k+1}-\min \eta_{k+1}}{\max\eta_{k+1}-\min \eta_{k+1}}$,
since $u$ is a characteristic function.
This choice is to control the rate of change of $\mathcal{A}[y(u_{k})]$. For the experiment, we employ one positive regional sources $f=x_1^2+0.1$, following \cite[Assumption 2.2]{beretta2018a}.

\begin{figure}[hbt!]
\centering
\begin{tabular}{ccc}
\includegraphics[width = \figlen, trim = {0.0cm 0.0cm 0.0cm 1.5cm}, clip]{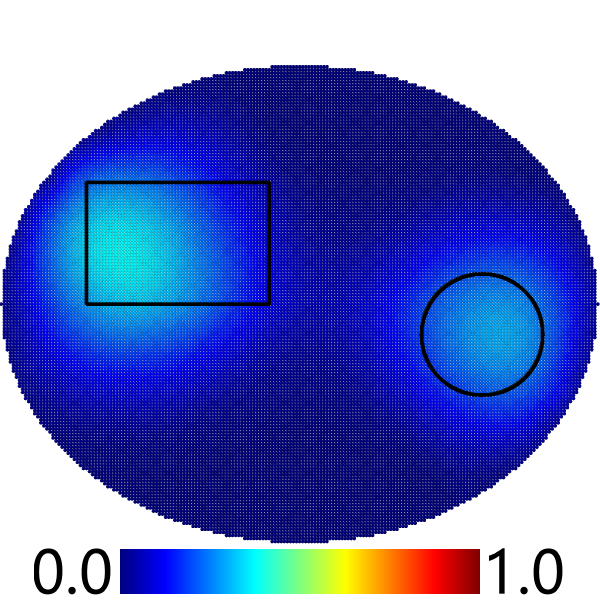}
&\includegraphics[width = \figlen, trim = {0.0cm 0.0cm 0.0cm 1.5cm}, clip]{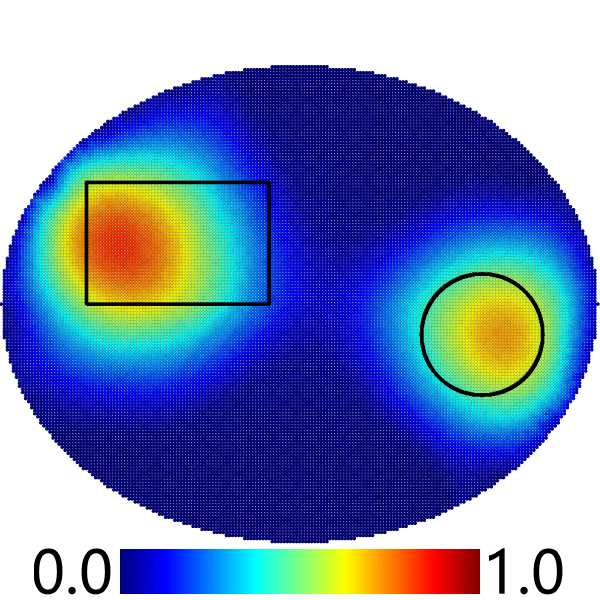}
&\includegraphics[width = \figlen, trim = {0.0cm 0.0cm 0.0cm 1.5cm}, clip]{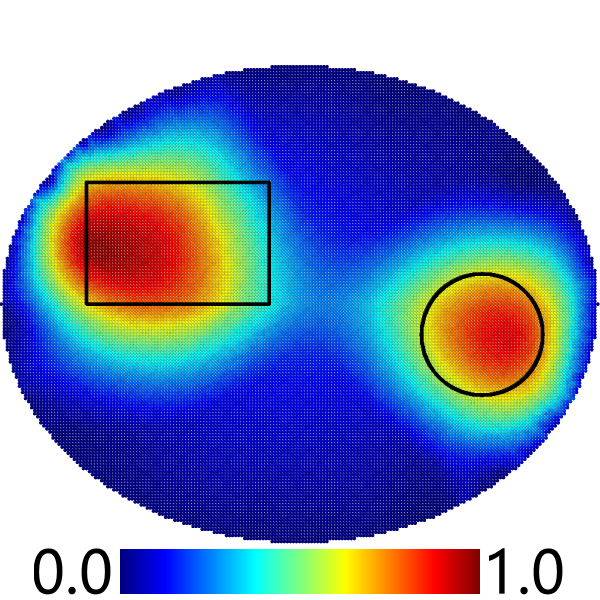}\\
\includegraphics[width = \figlen, trim = {0.0cm 0.0cm 0.0cm 1.5cm}, clip]{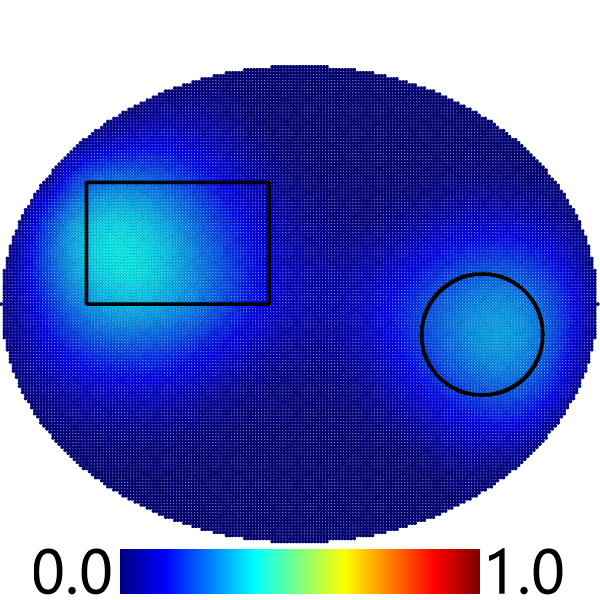}
&\includegraphics[width = \figlen, trim = {0.0cm 0.0cm 0.0cm 1.5cm}, clip]{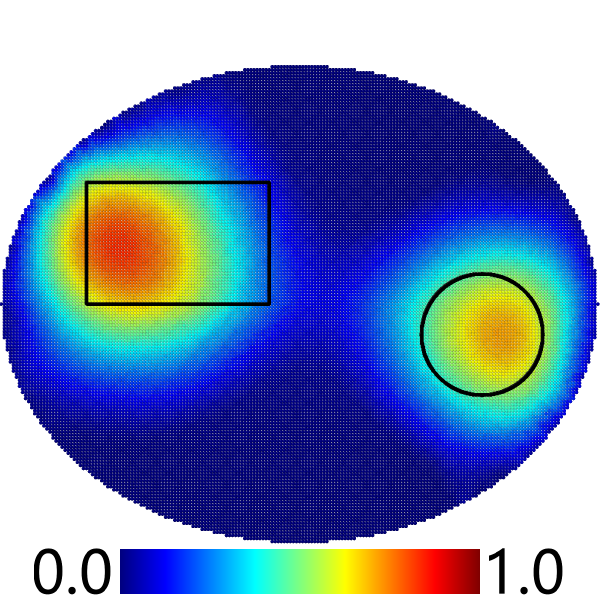}
&\includegraphics[width = \figlen, trim = {0.0cm 0.0cm 0.0cm 1.5cm}, clip]{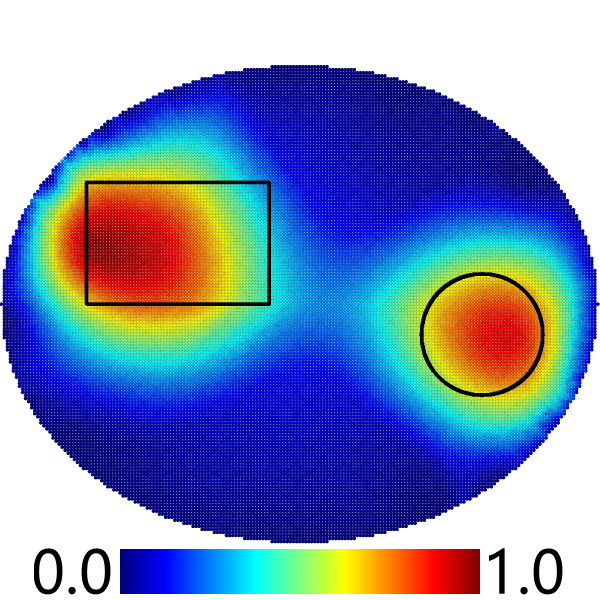}\\
(a) $k=1$ & (b) $k=6$ & (c) $k=11$
\end{tabular}
\caption{Visualization of $u_{k}$ for Example 4 with a noise level $\varepsilon = 10\%$.
The conductivity is $\text{1e-4}$ inside the inclusions and $1.0$ outside.
The top row displays the outcomes using the BFG correction, while the bottom row shows the results using the DFP correction.}
\label{fig4}
\end{figure}

The initial estimate $u_1$ (i.e., $k=1$) provides a baseline for the IDSM, which does not fully capture the inclusions; see Fig. \ref{fig4} for the reconstructions.
The IDSM continues to improve the accuracy of the inclusion estimate $u_k$:
the estimate $u_6$ (at the $6$th iteration) can already identify the smaller inclusion, and the estimate $u_{11}$ (at the $11$th iteration) further refines the recovery, giving a more precise and detailed localization.
The smaller inclusion that was initially overlooked is now well captured, highlighting the capability of the IDSM to enhance the resolution of the solution during the iteration.
These results underscore the effectiveness of the IDSM in handling nonlinear models, for which the iterative nature is especially beneficial.

\subsection{Example 5: non-smooth nonlinear model}
\begin{figure}[hbt!]
\centering
\begin{tabular}{ccc}
\includegraphics[width = \figlen, trim = {0.0cm 0.0cm 0.0cm 1.5cm}, clip]{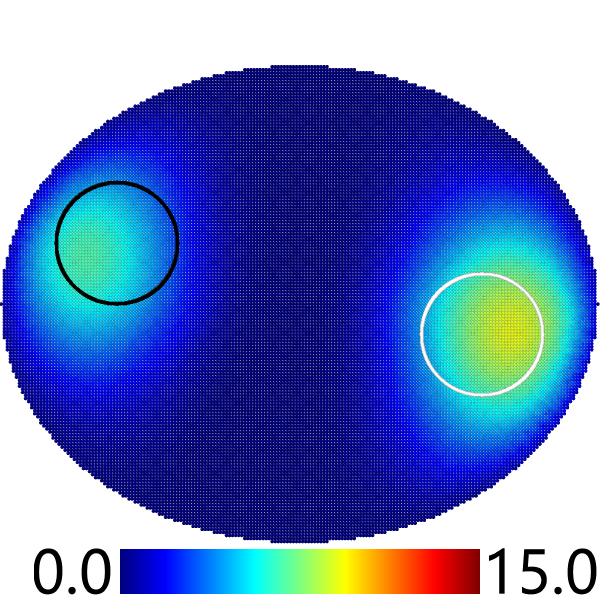}
&\includegraphics[width = \figlen, trim = {0.0cm 0.0cm 0.0cm 1.5cm}, clip]{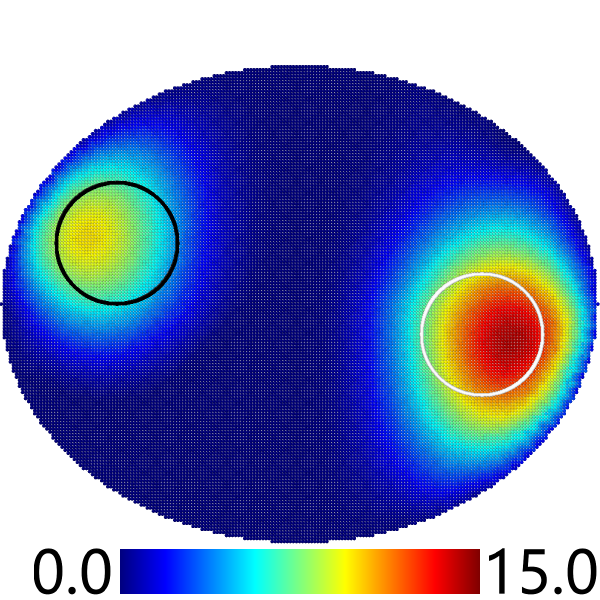}
&\includegraphics[width = \figlen, trim = {0.0cm 0.0cm 0.0cm 1.5cm}, clip]{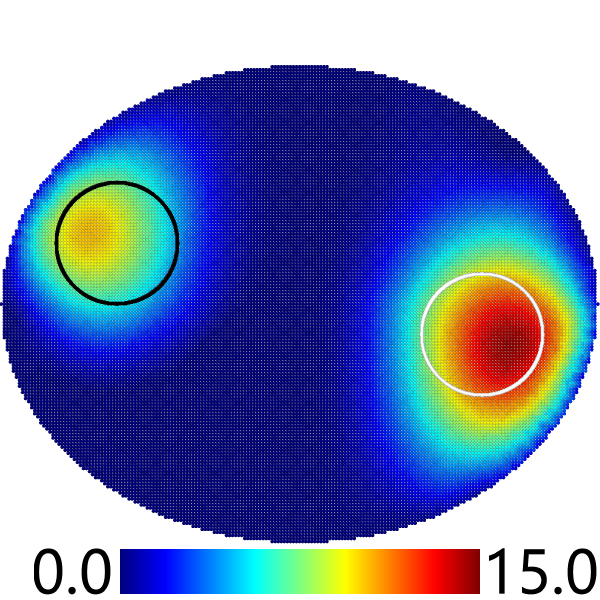}\\
\includegraphics[width = \figlen, trim = {0.0cm 0.0cm 0.0cm 1.5cm}, clip]{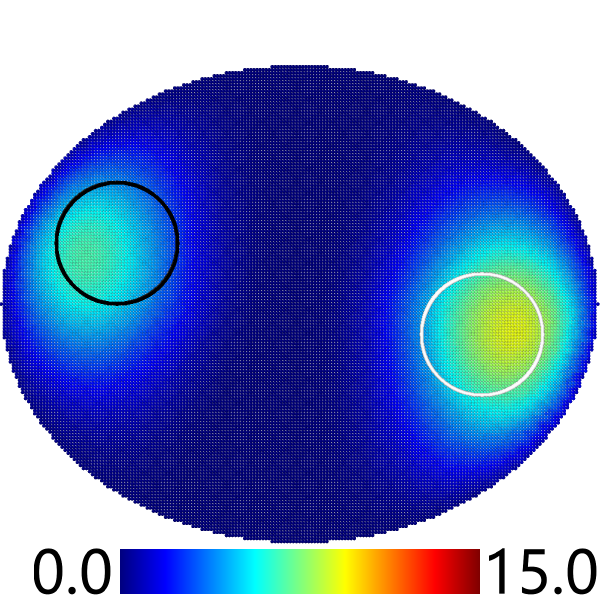}
&\includegraphics[width = \figlen, trim = {0.0cm 0.0cm 0.0cm 1.5cm}, clip]{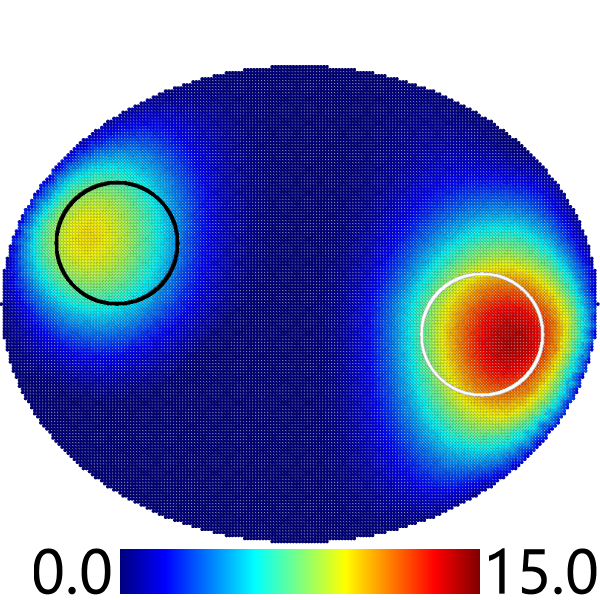}
&\includegraphics[width = \figlen, trim = {0.0cm 0.0cm 0.0cm 1.5cm}, clip]{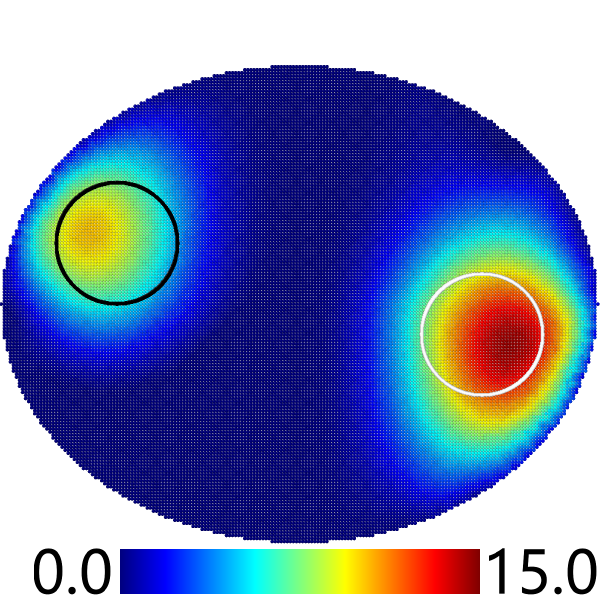}\\
(a) $k=1$ & (b) $k=6$ & (c) $k=11$
\end{tabular}
\caption{Visualization of $u_{k}$ for Example 5 with a noise level $\varepsilon = 10\%$.
The coefficient of $u$ is set to $40$ inside the white circle, $20$ inside the black circle, and $0$ outside both regions.
The top row displays the outcomes using the BFG correction, while the bottom row shows the results using the DFP correction.}
\label{fig5}
\end{figure}

In this example, we recover the potential in a non-smooth, nonlinear elliptic problem:
\begin{equation*}
\left\{
\begin{aligned}
-\Delta y + y + u|y|y &= f, \quad \text{in } \Omega, \\
\frac{\partial}{\partial n} y &= 0, \quad \text{on } \Gamma.
\end{aligned}
\right.
\end{equation*}
The associated operators $\mathcal{A}$ and $\mathcal{B}$ are defined, respectively, by
\begin{align*}
\langle \mathcal{A}y, w \rangle_{(H^{1}(\Omega))', H^{1}(\Omega)} &= \int_{\Omega} \nabla y \cdot \nabla w + yw\, \mathrm{d}x, \\
\langle \mathcal{B}[u]y, w \rangle_{(H^{1}(\Omega))', H^{1}(\Omega)} &= \int_{\Omega} u|y|yw\, \mathrm{d}x.
\end{align*}
Note that while the operator $\mathcal{B}[u]$ is nonlinear with respect to the state $y$, the operator $\mathcal{A}$ is actually linear in $y$. Thus Algorithm \ref{alg2} can be applied to the problem. In the experiment, we set the penalty parameter $\alpha=0.1$ in the regularized DtN map $\Lambda_{\alpha}(\mathcal{A})$ and utilize only one pair of Cauchy data corresponding to the regional source $f=x_{1}$.
We initialize the inhomogeneity $u$ to $u=0$ and the resolver $\mathcal{R}$ to $\mathcal{R}=d(x,\Gamma)$.
The map $\mathcal{P}$ is defined as a pointwise projection operator $\mathcal{P}(\eta)=\max(\min(\eta, 50.0), 0.0)$.
The results of the IDSM are shown in Fig. \ref{fig5}, where the parameter $u$ is $40$ and $20$ within the region enclosed by the white circle and black circle, respectively, and vanishes outside both regions.

Fig. \ref{fig5} presents the progression of the iterative process.
The initial estimate ($k=1$) provides a baseline for the IDSM, which is not sufficiently accurate in capturing the exact locations and values of the two inclusions.
By the $11$th iteration, the IDSM has effectively converged to a more accurate reconstruction.
The locations of both inclusions are now well-captured, with clear distinctions between the regions with $u=20$ and $u=40$.

\section{Concluding remarks}\label{sec_con}
By integrating an iterative mechanism into existing direct sampling methods, the proposed iterative direct sampling method  (IDSM) enjoys remarkable robustness against data noise and an enhanced capability to accurately identify inhomogeneities within a given medium.
It can iteratively refine the estimation of inclusion locations and shapes, leading to improved accuracy during the iteration process.
The numerical experiments on various models consistently show the superiority of the IDSM in handling noise and its flexibility in dealing with inhomogeneities of different types.
The framework is applicable to a wide array of elliptic problems, both linear and nonlinear, which makes the method a very promising tool for a broad spectrum of applications.
Future studies may focus on exploring its applicability to other types of inverse problems, e.g., geometric and boundary unknowns, and inverse problems associated with time-dependent PDEs, especially recovering moving inhomogeneities, which have posed significant difficulties to existing direct sampling methods.

\bibliographystyle{siam}
\bibliography{ref}
\end{document}